\documentclass[11pt]{article}
\usepackage{cite}
\usepackage[T2A]{fontenc}
\usepackage[utf8]{inputenc}
\usepackage[russian,english]{babel}

\usepackage{amsmath, amsfonts}

\usepackage{graphicx}
\usepackage{graphics}
\usepackage{amsmath}
\usepackage{amsthm}
\usepackage{amssymb}

\numberwithin{equation}{section}

\def\rd{{\rm d}} \def\re{{\rm e}} \def\rw{{\rm w}}

  \def\rF{{\rm F}}    
\def\rI{{\rm I}}  \def\rS{{\rm S}} 
\def\rT{{\rm T}} \def\rU{{\rm U}} 
\def\rGa{{\rm{Ga}}} \def\rdet{{\rm{det}}} \def\rGam{{\rm{Gam}}}
\def\rExp{{\rm{Exp}}} \def\rPo{{\rm{Po}}}
\def\rPDF{{\rm{PDF}}} \def\rRHS{{\rm{RHS}}}
\def\rTV{{\rm{TV}}} \def\rtr{{\rm{tr}}}

\def\bbE{\mathbb E} 
    \def\bbR{{\mathbb R}}

\def\bX{\mathbf X}  
\def\bx{\mathbf x}

\def\bb{\mathbf b}

 \def\bX{\mathbf X}

\def\bc{\mathbf c}

\def\cX{\mathcal X}

\def\FF{\mathfrak F}

\def\Gam{{\Gamma}} \def\Tha{{\Theta}}

\def\kap{\kappa}\def\lam{\lambda}
\def\gam{\gamma} \def\eps{\epsilon} 
 \def\vphi{{\varphi}}
 \def\vsig{{\varsigma}} 
\def\tha{{\theta}} \def\vtha{{\vartheta}} 
\def\ups{{\upsilon}}

\def\beq{\begin{equation}} \def\eeq{\end{equation}}
\def\beal{\begin{array}{l}} \def\beacl{\begin{array}{cl}}
\def\beac{\begin{array}{c}} \def\bear{\begin{array}{r}} 
\def\ena{\end{array}}
\def\onwl{\operatornamewithlimits} 

\def\diy{\displaystyle} 
\def\Blt{\blacktriangle}
 
\def\wvphi{{\widehat\vphi}} \def\obeta{{\overline\beta}}
\def\olam{{\overline\lambda}}  \def\omu{{\overline\mu}}
\def\otha{{\overline\theta}}  \def\osigma{{\overline\sigma}} \def\oSigma{{\overline\Sigma}}

\title{Context-sensitive hypothesis-testing and exponential families}
\author{Mark Kelbert\footnote{Honorary Reader in Statistics, University of Swansea, SA1 8EN, UK; The \\ 
National Research University ''The Higher School of Economics", Moscow 10100, RF} \ and 
Yuri Suhov\footnote{Math Dept, Penn State 
University, 16802 PA, US; University of Cambridge and St John's College, Cambridge, CB2 0WB and CB2 1TP, UK}}
\date{}

\begin{document}
\maketitle
\vskip .1truecm
{\bf Abstract.} We propose a number of concepts and properties related to `weighted' statistical inference
where the observed data are classified in accordance with a `value' of a sample string. The motivation 
comes from the concepts of weighted information and weighted entropy that proved useful in 
industrial/microeconomic and medical statistics. We focus on applications relevant in hypothesis testing and
an analysis of exponential families. Several notions, bounds and asymptotics are established, which generalize 
their counterparts well-known in standard statistical research. It includes Stein-Sanov theorem, Pinsker's, 
Bretangnole-Huber and van Trees inequalities and Kullback--Leibler, Bhattacharya, Bregman,
Burbea-Rao, Chernoff, R\'enyi and Tsallis divergences. 
\vskip .5cm 

\vskip .5cm

\section{Introduction}

This work focuses on a new aspect of `classical' statistics where equiprobable outcomes may have 
different `values' or `weights'.
Here one can speak of context-sensitive or context-dependent statistical observations and procedures 
such as estimation or hypothesis testing. An instructive example is provided by the Shannon entropy
and its derivatives and their `weighted' analogs. The Shannon entropy of a probability distribution with 
$m$ outcomes $x_1$, $\ldots$, $x_m$ and their probabilities $p(x_1)$, $\ldots$, $p(x_m)$ is given as
$h=-\sum\limits_ip(x_i)\ln p(x_i)$; here outcomes with the same probability provide equal contributions.
The weighted Shannon entropy can be defined as 
$$h^\rw_\vphi (p)=-\sum\limits_i\vphi (x_i)p(x_i)\ln p(x_i)$$ 
where $\vphi (x)$ represents the weight or value of outcome $x$ and $p=(p(x_i):1\leq i\leq m)$. 
A similar formulas can be employed in the case of a probability density:
\beq\label{eq1:wSE}h^\rw_\vphi (p)=-\int_\cX\vphi (x)p(x)\ln p(x)\nu (\rd x).\eeq
Here $\cX$ is a space of outcomes with a reference measure $\nu$ and $p$ is a probability 
density function, with $p(x)\geq 0$ and $\diy\int_\cX p(x)\nu (\rd x)=1$.

The significance of the weighted extension of the entropy can be illustrated by two examples
from industrial/microeconomic and clinical trial statistics. In papers \cite{FrSu1, FrSu2}
the authors considered a production line with a random input flow of raw or semi-fabricated
material. The usual/standard entropy of the input process measures variability of trajectories
of this process, and one might be tempted to predict that in a floor with many lines the bottleneck
of the global production process can be detected by comparing the empirical entropies for different
lines of the floor. However, a more accurate assessment is based on weighted entropies.
Here, high-quality items/units arriving at the input port of a production line relatively rarely are set 
to have a negligibly small weight since such items do not cause delays or disruptions
of the process. On the contrary, low-quality items/units that also occur relatively rarely, should be given a 
high value/weight (e.g., equal to $1$). This is motivated by the fact that low-quality input leads
to slowing the production down or even breaks of the whole line.
In the situation of complex industrial processes that are
mutually dependent, such a distinction helps to set correct priorities, particularly under limited resources
available for improving the quantity and quality of the output product.  A similar analysis was performed
in \cite{FrSu2} for industrial supply chains.  

Another example emerges in statistics of clinical trials; cf. \cite{JM, KKM1, KKM2}.
Suppose a doctor needs to choose between two drugs for a treatment
of a certain disease. The drugs are new and have not been tested as yet under the 
specific circumstances that take place in the trials under consideration. 
Assume that subsequent trials are independent. %(whatever it may mean
%in formal terms). 
In an imaginary case where the doctor is indifferent to possible results,
he/she may prescribes drug 1 and drug 2 intermittently. (Then the information about drugs' efficiency 
will be rapidly acquired.) 
Suppose that after a number of trials, it transpires that
drug 1 was considerably more effective than drug 2. Consequently, in a half of cases the doctor   
adopted a wrong decision. But we want to discover the potency of the drugs and simultaneously reduce
the fraction of wrong prescriptions.

Consequently, when dealing with trial $n+1$, the doctor must take into account results of previous $n$ 
trials. The goal is to achieve a given portion $\gamma$ of successful trials
(e.g., with $\gamma =0.95$.) The proposed mode of action  is as follows.  Suppose that after $n$ trials
there were $n_i\alpha_i(n_i)$ successes in $n_i$ attempts to prescribe drug $i$; here $i=1,2$ and
$n_1+n_2=n$. Maximizing the weighted Shannon entropy leads to the statistic
$$\sigma^{\rw ,\,\rS}_i=\dfrac{[\alpha_i(n_i)-\gamma ]^2}{2\alpha_i(n_i)[1-\alpha_i(n_i)]}(n_i+\beta +2)^{2\kap -1},
\;\;\kap\in [1/2,1],$$
while using a weighted Fisher information  
$$\sigma^{\rw ,\,\rF}_i=\dfrac{[\alpha_i(n_i)-\gamma ]^2}{\alpha_i(n_i)^2[1-\alpha_(n_i)]^2}
(n_i+\beta+2)^{2\kap},\;\;\kap\in [0,1].$$

The conclusion is to prescribe drug $i$ when $\sigma_i$ is smaller than the other; if they are 
equal toss the fair coin. Here $\beta >0$ and $\kap$ are parameters of the model.  
%A given patient may exhibit a considerable variability in observed 
%parameters of their life-functioning, which may seem alarming. However, an accurate assessment
%of the patient's state can be achieved by introducing appropriate quantitative weights to values
%of different parameters observed    

It turns out that the above procedure is asymptotically equivalent to minimization of the difference 
$h^\rw-h$ between a weighted entropy $h^\rw_\vphi (f)$ (see \eqref{eq1:wSE}) and its standard 
version $h=-\diy\int f(x)\ln\,f(x)\nu (\rd x)$ for recursively calculated weight function (WF) $\vphi =\vphi_{n_i}$ 
and probability density $f=f_{n_i}$. It was demonstrated that this algorithm gives better results 
versus other procedures known in the clinical trial practice in the sense of the simultaneous optimization 
the statistical power and the PCA (the probability of a correct allocation).

Successful applications of the notion of a weighted entropy motivate a study of mathematical
concepts or properties considered as weighted versions of well-known concepts and properties from
probability and statistics. In this paper we present a number of such instances in a rather simple form.
It may happen that further works will be done in this direction, extending existing limits of applicability 
of principles of statistical inference in response to growing demands, both in academic research 
and numerous applications.    

In fact, it is conceivable that a new, `weighted', or context-dependent, variant of statistical theory 
may emerge, where observed data have to be classified by a degree of importance, and 
procedures of statistical inference should take into account an additional `mark' indicating the
weight of a given sample.  For example, in testing a simple hypothesis versus a simple alternative
an error-probability should be modified to an (average) `cost' obtained by summing together
individual weights/values of sample strings which are included in a critical region. In this paper
we provide an initial result in this direction: see Lemma 2.2. Similarly, 
when we estimate a parameter from a sample, it makes sense to require that the estimator
is closer to the true value on samples with a high weight: this would lead to weighted distances
as a basis for measuring goodness of a given estimator.  Moreover, the approximation might 
be required to have better accuracy for some `anomalous' values of the parameter (e.g, a local 
temperature of an internal organ in a living body).

It is worthwhile to note that, from a formal point of view, introducing a weight function $\vphi$
means, simply speaking, that one passes from probability measures and densities to their
non-normalized counterparts. However, the point is to develop a theory that highlights the role 
of $\vphi$ in a maximally explicit form, tailor-made for applications.    

For this paper we selected a collection of concepts and properties with short and straightforward 
definitions and proofs. To make the presentation systematic, we focus on providing weighted 
analogs of (i) the Neyman--Pearson lemma assessing error-losses in hypothesis-testing, (ii) 
the Stein--Sanov lemma specifying log-asymptotics for type-II error-losses/costs, (iii) Cram\'er--Rao 
and van Trees inequalities for the mean-square deviation cost of an estimator
and (iv) explicit formulas for various weighted divergences 
in exponential families. The latter could be instructive for error-cost evaluations  arising in 
hypothesis testing. We did not pursue the goal of giving a complete picture in a single
paper and plan to continue the emerging line of research in future works.

In Section 2 we introduce weighted metrics and divergences, which are 
applied to the estimations of context-sensitive errors in hypothesis testing. This includes 
Neyman--Pearson type formulas 
for the minimal combined error-loss (Lemmas 2.2--2.4) and a weighted analog of the Stein--Sanov 
lemma (Theorem 2.1). Section 3 is dedicated to weighted Cram\'er--Rao and van Trees inequalities
which are presented in several versions (Theorems 3.1--3.2 and 3.3--3.5, respectively). 
In Section 4 a convenient formalism is explored, for computing weighted entropies and divergences 
in an exponential family. Conclusions from and discussion of our results, together with
possible future perspectives, are discussed in Section 5. 

Specific examples of calculations of weighted characteristics and parameters are given in Appendices
A--C to the main body of the paper. 

Throughout the paper, the term weighted and context-sensitive are used as equal in right.
Symbol $\Box$ marks the end of a proof and $\Blt$ the end of a definition.

\section{Context-sensitive hypothesis testing}

\hskip .9cm {\bf 2.1.} {\it Weighted total variation distance.} We consider probabilistic measures on an 
outcome space $\cX$ with a sigma-algebra $\FF$ and a given
(measurable) function $\vphi \geq 0$. For definiteness, it will be assumed that $(\cX,\FF)$ is a Polish space.
\vskip .5cm

{\bf Definition.}
Define the weighted total variation (TV) distance $\tau^\rw_{\vphi}(P,Q)$ under weight function (WF) $\vphi$,
between probability measures $P , Q $ on $(\cX,\FF),$ by  
\beq\label {10.0}\beacl\tau^\rw_{\vphi}(P,Q)&= \dfrac{1}{2}\left(\onwl{\sup}\limits_{A\in\FF}\left[\diy\int_A \vphi {\rd}P
 -\diy\int_A \vphi {\rd} Q \right]\right.\\ 
\;&\qquad\qquad\left.+\onwl{\sup}\limits_{A\in\FF}\left[\diy\int_A \vphi {\rd} Q -\diy\int_A \vphi 
{\rd}P \right]\right);\ena\eeq
сf. \cite{Kel}. When $\vphi \equiv 1$, this defines a standard TV distance.\quad$\Blt$\vskip .5cm

In what follows, we use an assumption that probability measures $P , Q $ are absolutely
continuous with respect to a measure $\nu$ and denote by $p$, $q$ their probability density functions (PDFs).
For definiteness, one may take $\nu =\nu (P,Q)=\dfrac{1}{2}\big(P+Q\big)$. Consequently, $\tau^\rw_\vphi (P,Q)$
will be referred to as $\tau^\rw_\vphi (p,q)$; the same convention will be applied to other 
functionals of the pair $P,Q$. Furthermore, arguments $p$, $q$ will be often omitted. We 
also write systematically $\diy\int f$ for an integral $\diy\int_\cX f\rd\nu$. If $f=\vphi g$, we adopt 
a convention to set 
\beq\label{eq:Evphi}\diy\int\vphi g=E_\vphi (g)\eeq
or even  $\diy\int\vphi g=E(g)$ when it is clear what WF $\vphi$ 
we address. The latter agreement will be applied to other functionals as well (e.g., in the notation 
$\tau^\rw (p,q)$). If $g=p$ is a PDF on $\cX$, $E_\vphi (p)$ yields the mean $\bbE_p[\vphi (X)]$ of the 
WF under $p$ (i.e., with $X\sim p$).
\vskip .5cm

{\bf Lemma 2.1.} (a) {\sl The following property holds true:}
\beq\beac\tau^\rw_{\vphi}(p,q)=\dfrac{1}{2}\diy\int\vphi\,|p-q|=\dfrac{1}{2}E_\vphi\big(|p-q|\big).
\ena\label {10.1}\eeq

(b) {\sl Suppose that $\vphi >0$ on $\cX$. Then $\tau^\rw_\vphi$ defines a metric for probability 
measures on $(\cX,\FF)$  i.e., for any 
{\rPDF}s $p$, $q$ and $r$:
{\rm{(i)}} $\tau^\rw_\vphi (p,q)\geq 0$ and $\tau^\rw_\vphi (p,q)=0$ iff $p$ and $q$ coincide $\nu$-a.s.
on $\cX$,  
{\rm{(ii)}} $\tau^\rw_\vphi (p,q)=\tau^\rw_\vphi (q,p)$, and {\rm{(iii)}} $\tau^\rw_\vphi (p,q)
\leq \tau^\rw_\vphi (p,r)+\tau^\rw_\vphi (r,q)$.}\vskip .5cm

{\it Proof.} (a) Consider the set $C^*=\big\{x\in\cX: q(x)\geq p(x)\big\}$ and write:
$$ E_\vphi \big((p-q){\bf 1}_{\cX\setminus C^*}\big)=\dfrac{1}{2}\big(E_\vphi (p)
-E_\vphi (q)\big)+\dfrac{1}{2}E_\vphi \big(|p-q|\big),$$
and
$$E_\vphi\big((q-p){\bf 1}_{C^*}\big) =\dfrac{1}{2}\big(E_\vphi (q)  
-E_\vphi (p)\big)+\dfrac{1}{2}E_\vphi \big(|p-q|\big).$$
%\eeq the set $C^*=\big\{x\in\cX: q(x)>p(x)\big\}$.
Summing these equalities one gets (\ref{10.1}). 

(b) Properties (i--iii) follow from representation \eqref{10.1}.
\qquad $\Box$ \vskip .5cm 

Some examples of calculation of the weighted TV distance are given in Appendix A.\vskip .5cm

{\bf 2.2.} {\it Context-sensitive error-losses in hypothesis testing.} In this section we relate the 
weighted TV distance 
to the sum of type I and II weighted error-losses in Neyman--Pearson hypothesis testing. 
Let $D:\;\cX\to [0,1]$ 
be a (possibly, randomized) decision rule for checking a simple hypothesis H$_0:\;p$ versus the 
simple alternative H$_1:\;q$. Define the type I and II error losses 
$\alpha^\rw_{\vphi}(D)$ and $\beta^\rw_{\vphi}(D)$:
\beq\label{albet:sc}\beal\alpha^\rw_{\vphi}(D)=\diy\int \vphi\,pD=E_\vphi (pD),\; 
\beta^\rw_{\vphi}(D)=\diy\int \vphi\,q(1-D)=E_\vphi (q(1-D)).\ena\eeq
Also, set: 
\beq\label{Delta:sc}
\beal\Delta_\vphi (p,q)=\dfrac{1}{2}\diy\int \vphi (p+q)=\dfrac{1}{2}
\big(E_\vphi (p)+E_\vphi (q)\big).\ena\eeq
\vskip .5cm

{\bf Lemma 2.2.} \label{WTVD} {\sl The following property is satisfied}
\beq \inf_{D}\big[\alpha^\rw_{\vphi}(D)+\beta^\rw_{\vphi}(D)\big]=\Delta_\vphi (p,q) -\tau^\rw_{\vphi}(p,q). 
\label{10}\eeq\vskip .5cm

{\it Proof.} As before, set $C^*=\big\{x\in\cX: q(x)>p(x)\big\}$. The result follows from the equality:
%$\forall$ $d\,:\;\cX\to [0,1]$  _\vphi 
\beq\beacl
E_\vphi (pD)+ E_\vphi (q(1-D))&=E_\vphi (q) +E_\vphi ((p-q)D) \\
\;&=E_\vphi (q)+E_\vphi\left(|p-q|D\Big({\bf 1}_{\cX\setminus C^*}-{\bf 1}_{C^*}\right)\Big).
\label{2.6}\ena\eeq
The RHS in \eqref{2.6} is minimized when $D={\bf 1}_{C^*}$. Consequently, 
$$\beal\onwl{\inf}\limits_{D}\Big[\alpha^\rw (D)+\beta^\rw(D)\Big]
=E_\vphi\left(q\right) -E_\vphi\big({\bf 1}_{C^*}(q-p)\big) 
=\Delta-\tau^\rw .\quad\Box \ena$$
\vskip 1cm

{\bf 2.3.} {\it Weighted divergences and distances between probability distributions.}
Weighted divergences and weighted distances between probability distributions
may be instrumental when one considers hypothesis testing for repeated observations. 
For definitions and properties of, and comments on standard (non-weighted) prototypes
of the concepts appearing in this and subsequent sections, cf. monographs \cite{GN, 
Pin,Tsy}  
and the bibliography therein. We also note original and recent papers \cite{Bha,Che,
Nie1,Nie2,NB, NN,Ren,Tsa,VEH}.

Let us begin with some basic definitions.
%Properties of tests on repeated observations require  
\vskip .5cm

{\bf Definition.} Given PDFs $p$ and $q$, 
define the weighted Hellinger distance $\eta^\rw_\vphi =\eta^\rw_{\vphi}(p,q)$ and
the weighted  Bhattacharyya coefficient (or weighted affinity) $\rho^\rw_\vphi =\rho^\rw_{\vphi}(p,q)$:
\beq\beac
\eta^\rw_{\vphi}:=\dfrac{1}{\sqrt{2}}\Big(E_\vphi\left((\sqrt{p}-\sqrt{q})^2\right)\Big)^{1/2},\\
\rho^\rw_{\vphi}:=E_\vphi\left(\sqrt{pq}\right)=\Delta_\vphi (p,q)-\eta^\rw_\vphi(p,q)^2.\ena\eeq
$\Blt$
\vskip .5cm

{\bf Lemma 2.3.} {\sl The following inequalities hold true:}
\beq\beacl 0&\leq
\Delta_\vphi (p,q)-\rho^\rw_{\vphi}(p,q)\\
\;&\leq \tau^\rw_{\vphi}(p,q)
\leq \sqrt{\Delta_\vphi (p,q)^2-\rho^\rw_{\vphi}(p,q)^2}.\ena\label{2.38}\eeq
\vskip .5cm

{\it Proof.} %For brevity, omit the arguments $p, q$. 
Clearly, $\rho^\rw_\vphi (p,q)\leq \Delta_\vphi (p,q)$.  
Further, the lower bound coincides with the lower bound in weighted  
Le Cam's inequality \cite{Kel}, Proposition 3.8. It can be proved as follows:
$$\beal
\Delta-\rho^\rw_\vphi  (p,q)=\dfrac{1}{2}E_\vphi\left((\sqrt{p}-\sqrt{q})^2\right)=\eta^\rw_\vphi (p,q)^2   \\
\qquad\leq\;\dfrac{1}{2}E_\vphi\big(|\sqrt{p}-\sqrt{q}|(\sqrt{p}+\sqrt{q})\big)   
=\dfrac{1}{2}E_\vphi (|p-q|)=\tau^\rw_\vphi  (p,q).\ena$$

The upper bound follows from the Cauchy-Schwarz inequality
$$\beal\dfrac{1}{4}E_\vphi (|p-q|)^2\leq \dfrac{1}{4}E_\vphi \big((\sqrt{p}
-\sqrt{q})^2\big) E_\vphi \big((\sqrt{p}+\sqrt{q})^2\big)\\
\quad = \dfrac{1}{4}\Big(2\Delta_\vphi (p,q)
-2E_\vphi (\sqrt{pq})\Big)\Big(2\Delta_\vphi (p,q)
+2E_\vphi (\sqrt{pq})\Big)\\
\quad =\Delta_\vphi (p,q)^2-(\rho^\rw (p,q))^2.\ena$$ \quad $\Box$ 
\vskip 1cm

Hence, in view of Lemma 2.2, we obtain\vskip .5cm

{\bf Lemma 2.4.} {\sl The following inequalities hold true:}
\beq\label{eq:009}\beal 0\leq\dfrac{\rho^\rw_\vphi (p,q)^2}{2\Delta_\vphi (p,q)}
\leq\Delta_\vphi (p,q)-\sqrt{\Delta_\vphi (p,q)^2-\rho^\rw_{\vphi}(p,q)^2}\\
\qquad\qquad\qquad\leq\onwl{\inf}\limits_{D}
\big[\alpha_{\vphi}(D)+\beta_{\vphi}(D)\big]\leq  \rho^\rw_{\vphi}(p,q).
\ena\eeq\vskip .5cm

{\it Proof.} For the lower bound, observe that for $a\geq b\geq 0$
$$a-\sqrt{a^2-b^2}\geq a-\sqrt{a^2-b^2+b^4/(4a^2)}$$
Setting $a=\Delta$, $b=\rho$ gives the lower bound. \quad $\Box$

\vskip 1cm

{\bf Definition.} Let $p$ and $q$ be two PDFs on $\cX$.
The weighted Kullback-Leibler divergence (of $p$ from $q$) is defined by  
\beq\label{eq:wKuLe}\beacl{\tt K}^\rw_{\vphi}(p\vert\vert q)
&=\diy\int \vphi\,p{\bf 1}(p>0)\,\ln\big(p/q\big)=E\Big(p{\bf 1}(p>0)\ln\big(p/q\big)\Big)\\
\;&=\bbE_p\Big[\vphi (X){\bf 1}(p(X)>0)\ln\big(p(X)/q(X)\big)\Big].\ena\eeq
Here and below, $\bbE_p$ denotes expectation 
under PDF $p$. \quad $\Blt$\vskip .5cm

{\bf Lemma 2.5.} {\sl The condition $E(p)\geq E(q)$ implies the weighted 
Gibbs inequality 
\beq{\tt K}^\rw_{\vphi}(p\vert\vert q)\geq 0,\eeq
with equality iff $\vphi p=\vphi q$ \ $\nu$-a.s.} \vskip .5cm
 
{\it Proof.} See [KSSY], Theorem 1.3.\quad$\Box$\vskip .5cm

The weighted Kullback--Leibler divergence allows us to lower-bound the RHS in Lemma 2.2:
\vskip .5 cm

{\bf Lemma 2.6.} (The weighted Pinsker  inequality)  {\sl Suppose $E_\vphi (p)\geq E_\vphi (q)$. Then
the weighted {\rTV} distance obeys}
\beq\tau^\rw_\vphi (p,q)\leq \sqrt{{\tt K}^\rw_{\varphi}(p||q)/2\,}\,\sqrt{E_\vphi (p)}.\label{10.3}\eeq

{\it Proof.} Cf. \cite{Tsy}, Lemma 2.5. Define the function $\Pi (x)=x\ln x-x+1$ for $x>0$. The 
following bound is claimed:  
\beq \Pi (x)=x\ln x-x+1\geq \dfrac{3}{2}\dfrac{(x-1)^2}{x+2}, \quad x>0.
\label{10.4}\eeq
Indeed, at $x=1$ both sides in \eqref{10.4} coincide, and their first derivatives coincide as well. 
Then the bound $\Pi''(x)=\dfrac{1}{x}\geq \dfrac{27}{(x+2)^3}$ proves the claim.

Now, assume that the supports of functions $\vphi p$ and $\vphi q$ coincide (otherwise the right-hand 
side in \eqref{10.3} is infinite, and the bound holds trivially). By the Cauchy-Schwarz inequality
$$\beal (\tau^\rw_\vphi (p,q))^2=\dfrac{1}{4}E_\vphi (q(p/q-1))^2\\
\quad\leq\dfrac{1}{4} E_\vphi\left(q\dfrac{(p/q-1)^2}{p/q+2}\right)E_\vphi\big(q(p/q+2)\big)\\
\quad\leq\dfrac{3}{4}E_\vphi\left(q\dfrac{(p/q-1)^2}{p/q+2}\right) E_\vphi (p)\\
\quad\leq\dfrac{1}{2}E_\vphi\big(q\Pi (p/q)\big)E_\vphi (p)\leq\dfrac{1}{2}E_\vphi (p)
{\tt K}^\rw_\vphi (p||q) . \quad \Box
\end{array}$$

%It is instrumental to note that the connection between the Fisher information and a `curvature' of a family
%of PDFs is preserved in the weighted case.  
\vskip .5 cm

{\bf Definition.} Given $0<\alpha <1$, the quantity 
\beq\label{ChernCo}\beal C^\rw_{\vphi ,\alpha}(p,q)=\dfrac{E_\vphi\left(p^{\alpha}q^{1-\alpha}\right)}{E_\vphi (p)}
\ena\eeq 
is called the weighted Chernoff $\alpha$-coefficient (between PDFs $p$ and $q$). 
The weighted Chernoff, R\'enyi and Tsallis $\alpha$-divergences ${\tt C}^\rw_{\vphi ,\alpha}(p,q)$, 
${\tt R}^\rw_{\vphi ,\alpha}(p,q)$, and ${\tt T}^\rw_{\vphi ,\alpha}(p,q)$
are defined by
\beq\label{1.3.11}\beac {\tt C}^\rw_{\vphi ,\alpha}(p,q)=-\ln C^\rw_{\vphi ,\alpha}(p,q),\;
{\tt R}^\rw_{\vphi ,\alpha}(p,q)=\dfrac{E_\vphi (p)}{1-\alpha}\ln C^\rw_{\vphi ,\alpha}(p, q),\\
{\tt T}^\rw_{\vphi ,\alpha}(p,q)=\dfrac{E_\vphi (p)}{1-\alpha}\Big(C^\rw_{\vphi ,\alpha}(p,q)-1\Big).
\ena\eeq
For brevity, we often omit the term weighted and the reference to the value of $\alpha$. Note that 
as $\alpha \to 1$, 
$${\tt R}^\rw_{\vphi ,\alpha}(p,q)\to {\tt K}^\rw_{\vphi}(p\| q),\;\;
{\tt T}^\rw_{\vphi ,\alpha}(p,q)\to {\tt K}^\rw_{\vphi}(p\vert\vert q).$$

A special role is played by the value $\alpha =1/2$. In particular, $C^\rw_{\vphi ,1/2}(p, q)
=\dfrac{1}{E_\vphi (p)}\rho^\rw_\vphi (p,q)$ whereas the quantity ${\tt C}^\rw_{\vphi, 1/2}(p,q)$ gives the 
Bhattacharyya divergence ${\tt A}^\rw_{\vphi}(p,q)$:
\beq\beacl{\tt A}^\rw_{\vphi}(p,q)&:={\tt C}^\rw_{\vphi, 1/2}(p,q) 
=-\ln\,\Big[\Delta_\vphi (p,q)-\eta^\rw_{\vphi}(p,q)^2\Big]+\ln E_\vphi (p )\\
\;&\;=-\ln\,\rho^\rw_\vphi (p,q)
+\ln E_\vphi (p ).\ena\eeq

The weighted Shannon entropy (see \eqref{eq1:wSE}) is defined as
\beq\label{eq2:wSE}
h^\rw_\vphi (p)=-\int\vphi p\ln p=E_\vphi (-p\ln p)
\eeq
while the weighted R\'enyi $\alpha$-entropy $r^\rw_{\vphi ,\alpha}(p )$ by
\beq\label{RenyiEnt}
r^\rw_{\vphi,\alpha}(p)=\dfrac{E_\vphi (p)}{1-\alpha}
\ln\dfrac{E_\vphi (p^{\alpha})}{E_\vphi (p)},
\quad 0<\alpha <1.\eeq
As in a standard (non-weighted) case,  
$r^\rw_{\vphi ,\alpha}(p)$ is expected to converge, as $\alpha\to 1$, to 
$h^\rw_{\vphi}(p)$. This can be verified 
under some regularity conditions upon $\vphi$ and $p$, when  
the limit can be carried under the integration in $E_\vphi (p^\alpha )$. We comment on this property
below, when we consider specific situations and examples. 
\vskip .5cm

An extended version of the weighted R\'enyi entropy can also be considered:
\beq\label{RenyiEntExt}\beacl
r^\rw_{\vphi,\alpha,\beta}(p)&=\dfrac{E_\vphi (p)}{1-\alpha}
\ln\dfrac{E_\vphi (p^{\alpha+\beta-1})}{E_\vphi (p^{\beta})},\\
\;&\quad \alpha >0, \alpha\neq 1,\; \beta>0, \alpha+\beta>1.\ena\eeq
%In this paper we do not discuss the Tsallis or
%extended R\'enyi entropies leaving them for separate publications. 
\quad $\Blt$ 
\vskip .5cm

{\bf 2.4.} {\it The weighted Bretagnolle–Huber bound.} The Bretagnolle–Huber
bound \cite{BH} strengthens the Pinsker inequality for the case where the Kullback-Leibler divergence 
is large.\vskip .5cm

{\bf Lemma 2.7.} (The Bretagnolle–Huber bound) {\sl The following inequality holds true:}
\beq \tau^\rw_{\vphi}(p,q)\leq \sqrt{\Delta^\rw_\vphi (p,q)^2-\re^{-{\tt K}^\rw_{\vphi}(p\|q)}}.
\label{17}\eeq
\vskip .5cm   

{\it Proof.} In view of \eqref{10.1},
$$\beac\tau^\rw_\vphi (p,q) =\Delta^\rw_\vphi (p,q)
-E_\vphi \big(\min[p,q]\big)=E_\vphi \big(\max[p,q]\big) -\Delta^\rw_\vphi (p,q) .\ena$$
Hence, by Cauchy-Schwarz and Jensen's inequalities,
$$\beal
\Delta^\rw_\vphi (p,q)^2-\tau^\rw_\vphi (p,q )^2
=E_\vphi \big(\min[p,q]\big)E\big(\max[p,q]\big)\\
\qquad\geq\Big(E_\vphi \Big(\sqrt{ \min[p,q]\max[p,q]}\Big)\Big)^2
=\Big(E_\vphi \big(\sqrt{pq}\big)\Big)^2\\
\qquad =\re^{2\ln E_\vphi (p\sqrt{q/p})}\geq \re^{2E_\vphi (p\ln\sqrt{q/p})}
=\re^{-{\tt K}^{\rw}_\vphi{(p||q)}}.\qquad \square\ena$$
\vskip .5cm

{\bf Remark.} Bound \eqref{17} improves on the weighted Pinsker inequality 
when ${\tt K}^\rw_\vphi (p\|q)$ is large.  
\vskip .5cm

{\bf 2.5.} {\it Weighted tests for subsequent observations.}
If the observer believes that s/he is dealing with IID observations $\bX_1^n=(X_1,\ldots ,X_n)$, 
it is natural to suppose that the WF $\bx_1^n\mapsto\vphi^{(n)}(\bx_1^n)$ factorizes:  
$\vphi^{(n)}(\bx_1^n)= \prod\limits_{i=1}^n \vphi (x_i)$, for $\bx_1^n=(x_1,\ldots ,x_n)\in\cX^n$.
We will follow the assumptions of independence and factorization throughout this subsection 
without stressing it every time again. More precisely, let $p^{(n)}(\bx_1^n)=\prod\limits_{i=1}^np(x_i)$
and $q(\bx_1^n)=\prod\limits_{i=1}^nq(x_i)$ be two PDFs relative to measure $\nu^n$ 
on $(\cX^n,\FF^n)$. As above, adopt the convention of
writing $\diy\int f$ instead of $\diy\int_{\cX^n} f\rd\nu^n$; also put 
\beq\label{eq:Envphi}E_{\vphi^{(n)}}(g)=\diy\int\vphi^{(n)}g\rd\nu^n\;\hbox{ and }\;\bbE_{p^{(n)}}[\vphi^{(n)}g]
=\diy\int\vphi^{(n)}gp^{(n)}\rd\nu^n\eeq 
for a given function $g:\;\cX^n\to\bbR$.  Consider a decision rule 
$D:\;\cX^n\to [0,1]$
for testing a simple hypothesis H$_0:\;p^{(n)}$ versus a simple alternative H$_1:\;q^{(n)}$\,. 
Following \eqref{albet:sc}, define the 
corresponding type I and II  weighted error-losses under rule $D_n$:
\beq\beac\alpha_{\vphi^{(n)}}(D_n)=E_{\vphi^{(n)}}(p^{(n)}D_n) ,\quad  
\beta_{\vphi^{(n)}}(D_n)=E_{\vphi^{(n)}}(q^{(n)}(1-D_n))\,.\ena\eeq
\vskip.5cm

The following lemmas give bounds for the minimal
sum $\alpha^\rw_{\vphi^{(n)}}(D_n)+\beta^\rw_{\vphi^{(n)}}(D_n)$. The lower and upper bounds 
in Lemma 2.8 hold for 
all $n$. The lower bound is asymptotically improved in Lemma 2.9.\vskip .5cm
 
{\bf Lemma 2.8.} {\sl The following bounds hold true:}
\beq\beal\dfrac{\rho^\rw_\vphi (p,q)^{2n}}{E_\vphi (p)^n+E_\vphi (q)^n}
\leq\onwl{\inf}\limits_{D_n}\Big[\alpha^\rw_{\vphi^{(n)}}(D_n)+\beta^\rw_{\vphi^{(n)}}(D_n)\Big]
\leq \rho^\rw_{\vphi}(p ,q)^n\\
\qquad\qquad =E_\vphi (p)^n\re^{-n{\tt A}^\rw_{\vphi}(p,q)}\leq \big(\Delta^\rw_\vphi (p,q)^2-\tau^\rw_{\vphi}(p,q)^2\big)^{n/2}.\ena
\label{2.17}\eeq
{\sl Moreover, if $\Delta^\rw_\vphi (p,q)\leq 1$ then}
\beq
\onwl{\inf}\limits_{D_n}\Big[\alpha^\rw_{\vphi^{(n)}}(D_n)+\beta^\rw_{\vphi^{(n)}}(D_n)\Big]
\leq  \re^{-n\eta^\rw_\vphi (p,q)^2}.\label{3.24}\eeq\vskip .5cm

{\it Proof.} The bounds \eqref{2.17} follow from \eqref{eq:009}. It remains to check \eqref{3.24}. Here we 
employ an elementary bound:
$$\ln \rho^\rw_{\vphi}(p,q)\leq \rho^\rw_{\vphi}(p,q)-1
=\Delta^\rw_\vphi (p,q)-1-\eta^\rw_\vphi (p,q)^2.
\label{3.25}$$

Under condition $\Delta^\rw_\vphi (p,q)\leq 1$, we obtain
$(\rho^\rw_\vphi (p,q))^n \leq\re^{-n \eta^\rw_{\vphi}(p ,q)^2}$, and \eqref{3.24} follows. \quad $\Box$
\vskip.5cm

{\bf Lemma 2.9.} {\sl Assume that $E_\vphi (p)\geq 1$. Then for $\forall \epsilon>0$ 
for $n>n_0(\eps)$
\beq\onwl{\inf}\limits_{D_n}\Big[\alpha_{\vphi^{(n)}}(D_n)+\beta_{\vphi^{(n)}}(D)\Big]
\geq \dfrac{1-\epsilon}{2\Delta^\rw_\vphi (p,q)}
\re^{-nE_\vphi (p)^{n-1}{\tt K}^\rw_{\vphi}(p ||q )}.\label{2.25}\eeq}
\vskip .5cm

{\it Proof.} In view of \eqref{10} and (\ref{17}), it seems reasonable to 
use a factorization of Kullback-Leibler divergence  %{\bf ugly!!}
\beq E_{\vphi^{(n)}}\left(p^{(n)}(\bx)\left(\sum\limits_{i=1}^n{\mathbf 1}(q(x_i)>0)
\ln\dfrac{p(x_i)}{q(x_i)}\right)\right)=nE_\vphi (p)^{n-1}{\tt K}^\rw_{\vphi}(p||q).\eeq
For  $n>n_0(\eps)$ we obtain by the Taylor expansion:
$$\beal\onwl{\inf}\limits_{D_n}\Big[\alpha_{\vphi^{(n)}}(D_n)+\beta_{\vphi^{(n)}}(D_n)\Big]
=\Delta^\rw_\vphi (p,q)-\tau^\rw_{\vphi^{(n)}}(p^{(n)},q^{(n)})\\
\quad \geq \Delta^\rw_\vphi (p,q)-\sqrt{\Delta^\rw_\vphi (p,q)^2
-\re^{-nE_\vphi (p)^{n-1}{\tt K}^\rw_\vphi (p||q)}}\\
\quad\geq \dfrac{1-\epsilon}{2\Delta^\rw_\vphi (p,q)}\re^ {-n E_\vphi (p)^{n-1}{\tt K}^\rw_{\vphi}(p||q)}.
\quad\square\ena$$ 
\vskip .5cm

{\bf 2.6.} {\it Weighted Stein--Sanov lemma/theorem.} Let $p,q$ be PDFs on $\cX$. As above, we want to 
test the simple hypothesis H$_0: p$ against the simple alternative  H$_1: q$ for an 
IID sample  $\bX_1^n=(X_1,\ldots ,X_n)\in {\cal X}^n$ and consider a decision rule
$D_n:\;\cX^n\to [0,1]$. Given  $\alpha\in (0,1)$, define 
\beq\beal\vsig_{\vphi^{(n)}}(\alpha )=\onwl{\inf}\limits_{D_n}\left[E^{(n)}\big((1-D)q^{(n)}\big)\right];
\quad\inf\;\hbox{over rules}\\ 
D_n:\;\cX^n\to [0,1]\;
\hbox{ with }\diy\int\vphi^{(n)}p^{(n)}D_n\leq E^{(n)}\alpha.\ena\eeq

The following theorem extends the result attributed to Sanov \cite{San} (also called
the Chernoff--Stein lemma) to the context-dependent case. \vskip .5cm  

{\bf Theorem 2.1.} {\sl Suppose that ${\tt K}^\rw_\vphi (p\|q) <\infty$. Then for any 
$\alpha\in (0,1)$
\beq\label{eq:wStSaY}\beal
\limsup\limits_{n\to\infty}\dfrac{1}{n}\ln \vsig_{\vphi^{(n)}}(\alpha )\leq
\ln\,E_\vphi (p)-{\tt K}^\rw_{\vphi}(p\|q )/E_\vphi (p),\\
\liminf\limits_{n\to\infty}\dfrac{1}{n}\ln \vsig_{\vphi^{(n)}}(\alpha )\geq
\ln\,E_\vphi (p)-{\tt K}^\rw_{\vphi}(p\|q )/E_\vphi (p)\ena\eeq
and hence,
\beq \lim\limits_{n\to\infty}\dfrac{1}{n}\ln \vsig_{\vphi^{(n)}}(\alpha )=
\ln\,E_\vphi (p)-{\tt K}^\rw_{\vphi}(p\|q )/E_\vphi (p).\eeq}

{\it Proof.} 
We want to upper- and lower-bound the integral $\diy\int_{\cX^n}\vphi^{(n)}q^{(n)}D_n$ 
for a suitable decision rule $D_n$. Set: 
$$\pi :=\vphi\,p/E_\vphi (p),\quad\vtha:=\vphi\,q/E_\vphi (q):\quad\hbox{ tilted one-time PDFs,}$$
$$\pi^{(n)}:=\vphi^{(n)}p^{(n)}/E_\vphi (p)^n,\;\vtha^{(n)}:=\vphi^{(n)}q^{(n)}/E_\vphi (q)^n:
\;\hbox{ tilted $n$-fold  PDFs.}$$ 
Next,
$$\beac z_i(\bx_1^n)=\ln\,\big[p(x_i)/q(x_i)\big],\;\;1\leq i\leq n,\\
\hbox{with}\;\bbE_{\pi^{(n)}}[z_i(\bX_1^n)]=\dfrac{{\tt K}^\rw_\vphi (p\|q)}{E_\vphi (p)},\;\;
\bbE_{\vtha^{(n)}}[z_i(\bX_1^n)]=-\dfrac{{\tt K}^\rw_\vphi (q\|p)}{E_\vphi (q)}.\ena$$ 
%and it is assumed that ${\tt K}^\rw_\vphi (p\|q)/E(p)\neq -{\tt K}^\rw_\vphi (q\|p)/E(q)$.
The key point is to use the LLNs for 
$z_i(\bx_1^n)$ under PDF $\pi^{(n)}$. 
Write:
\beq\beal\diy\int\vphi^{(n)}q^{(n)}D_n=\int\vphi^{(n)}p^{(n)}D_n\,\dfrac{q^{(n)}}{p^{(n)}}\\
\quad ={\diy\int}\vphi^{(n)}p^{(n)}D_n\re^{-\ln\,p^{(n)}/q^{(n)}}
=E_\vphi (p)^n{\diy\int}\pi^{(n)}D_n\exp\Big(-\sum\limits_{i=1}^nz_i\Big).\ena\eeq

By the LLNs, as $n\to\infty$, 
\beq\beac\hbox{under $\pi^{(n)}$}:\quad\dfrac{1}{n}\sum\limits_{i=1}^nz_i(\bX_1^n)
\to\dfrac{{\tt K}^\rw_\vphi (p\|q)}{E_\vphi (p)},\\
\hbox{under $\vtha^{(n)}$}:\quad\dfrac{1}{n}\sum\limits_{i=1}^nz_i(\bX_1^n)
\to-\dfrac{{\tt K}^\rw_\vphi (q\|p)}{E_\vphi (q)}\ena\eeq
Convergence in probability under $\pi^{(n)}$ means that $\forall$ $\eta >0$, 
$\diy\int\pi^{(n)}{\bf 1}_{L_n(\eta )}\to 1$. Here 
\beq\beac L_n (\eta ):=\left\{\bx_1^n\in\cX^n:\;-\eta\leq\dfrac{1}{n}\sum\limits_{i=1}^nz_i(\bx_1^n)
-\dfrac{{\tt K}^\rw_\vphi (p\|q)}{E_\vphi (p)}\leq\eta\right\}.\ena\eeq

Thus, if we take $D_n(\eta):=1-{\bf 1}_{L_n(\eta )}$, we obtain a decision rule with $\diy\int\vphi^{(n)}
p^{(n)}D_n\leq\alpha E_\vphi (p)^n$. For this rule: 
\beq\beacl{\diy\int}q^{(n)}(1-D_n)
&=E_\vphi (p)^n{\diy\int}\pi^{(n)}{\bf 1}_{L_n(\eta )}\exp\,\Big(-\sum\limits_{i=1}^nz_i\Big)\\
%\;&=E_\vphi (p)^n{\diy\int}\pi^{(n)}{\bf 1}_{L_n(\eta )}\exp\,\Big(-\sum\limits_{i=1}^nz_i\Big)\\
\;&\leq E_\vphi (p)^n\re^{n(-{\tt K}^\rw_\vphi (p\|q)/E_\vphi (p)+\eta)}\diy\int \pi^{(n)}
{\bf 1}_{L_n(\eta)}.\ena\eeq
This implies that 
$$\beac\limsup\limits_{n\to\infty}\dfrac{1}{n}\ln \vsig_{\vphi^{(n)}}(\alpha )
\leq\limsup\limits_{n\to\infty}\dfrac{1}{n}\ln\,{\diy\int}\vphi^{(n)}q^{(n)}{\bf 1}_{L_n(\eta )}\\
\qquad\leq \ln\,E_\vphi (p)-\dfrac{{\tt K}^\rw_\vphi (p\|q)}{E_\vphi (p)}+\eta\ena$$
which leads to upper bound in \eqref{eq:wStSaY}.

For the lower bound: assume that rule $D_n$ has 
$$\bbE_{\pi^{(n)}}[1-D_n]=\diy\int (1-D_n)\pi^{(n)}\geq 1-3\alpha $$
where $0\leq\alpha <1/3$. 
Then take $A_n:=\{\bx_1^n:\;D_n(\bx_1^n)\leq\alpha/3\}$ and deduce that 
$$\bbE_{\pi^{(n)}}[{\bf 1}_{A_n}]=\int{\bf 1}_{A_n}\pi^{(n)}\geq\dfrac{2 }{3}.$$ 
Write:
$$\beal\diy\int\vphi^{(n)}q^{(n)}(1-D_n)
\geq\left(1-3\alpha\right)\int\vphi^{(n)}q^{(n)}{\bf 1}_{A_n}{\bf 1}_{L_n(\eta )}\\
\qquad =\left(1-3\alpha\right){\diy\int}\vphi^{(n)}p^{(n)}{\bf 1}_{A_n\cap
L_n(\eta )}\exp\Big(-\sum\limits_{i=1}^nz_i\Big)\\
\qquad\geq\left(1-3\alpha\right)E_\vphi (p)^n\re^{n(-{\tt K}^\rw_\vphi (p\|q)/E_\vphi (p)-\eta)}
{\diy\int}\vphi^{(n)}p^{(n)}{\bf 1}_{A_n\cap L_n(\eta )}\\
\qquad =\left(1-3\alpha\right)E_\vphi (p)^n\re^{n(-{\tt K}^\rw_\vphi (p\|q)/E_\vphi (p)-\eta)}
{\diy\int}\pi^{(n)}{\bf 1}_{A_n\cap L_n(\eta )}.\ena$$
For $n$ large enough,
$${\diy\int}\pi^{(n)}{\bf 1}_{A_n\cap L_n(\eta )}\geq\dfrac{1}{3}.$$
Taking the log and passing to the limit leads to the lower bound in \eqref{eq:wStSaY}. 
\quad $\Box$
\vskip .5cm

\section{Weighted Cram\'er-Rao and van Trees inequalities}\vskip .5cm

\hskip .7cm In this section we consider a family of PDFs $p_\tha$ depending smoothly on 
a parameter $\tha =(\tha_1,\ldots ,\tha_d)\in\Theta$ where $\Theta\subseteq\bbR^d$ is a 
$d$-dimensional domain. In some places we adopt a simplified version, taking $d=1$,
to shorten calculations and make the formulas more transparent. Formally, the property of 
smoothness 
means that $\forall\;\tha\in\Tha$ there exists an $\FF$-measurable gradient function  
$$x\in\cX\mapsto\nabla_\tha p_\tha (x)=\big((p_\tha)'_{\tha_1}(x),\ldots , (p_\tha)'_{\tha_d}(x)\big).$$ 
%where $\nabla$ stands for the row-gradient $\nabla_\tha$; 
Here and below $(p_\tha)'_{\tha_l}$ stands for the partial derivative $\partial_{\tha_l}p_\tha$, 
$1\leq l\leq d$. The subscript 
$\tha$ in $\nabla_\tha$ will be omitted.  For $d=1$ we use $(p_\tha)'$, the derivative in $\tha$. 
For simplicity, it will be assumed that the values $\nabla p_\tha (x)$
and $(p_\tha)'(x)$ are defined $\forall$ $x\in\cX$. % and $x\mapsto\nabla p_\tha (x)$ and 
%$x\mapsto (p_\tha )'(x)$ are measurable functions $\forall$ $\tha\in\Tha$. 
We will need certain
assumptions upon $\nabla p_\tha$ and $(p_\tha)'$ which are specified in individual assertions below. 
\vskip .5cm 

{\bf Definition.}
The weighted Fisher information matrix is defined by
\beq\label{eq:wFiM}\beacl
\rI^\rw_\vphi (\tha )&=\diy\int\vphi\,{\mathbf 1}(p_\tha >0)
p_{\tha}^{-1}\nabla p_\tha^{\rT}\nabla p_\tha \\
\;&=\bbE_\tha\Big[\vphi (X){\mathbf 1}(p_\tha (X)>0)p_\tha (X)^{-2}\nabla p_\tha (X)^{\rT}
\nabla p_\tha (X)\Big].\ena\eeq 
Here, and below $\bbE_\tha$ denotes expectation 
under $p_\tha$. %Accordingly, a random variable $X$ under expectation $\bbE_\tha$
%has the PDF $p_\tha$. 
When $\vphi\equiv 1$, we obtain its non-weighted prototype
\beq\label{eq:sFIM}\beacl\rI (\tha )&=\diy\int\vphi\,{\mathbf 1}(p_\tha >0)
p_{\tha}^{-1}\nabla p_\tha^{\rT}\nabla p_\tha \\ 
\;&=\bbE_\tha\Big[{\mathbf 1}(p_\tha (X)>0)p_\tha (X)^{-2} 
\nabla p_\tha (X)^{\rT}\nabla p_\tha (X)\Big].\ena\eeq  

In the case $d=1$, we have scalar 
%parameter $\tha$, and $\Theta$ is an interval on
%a line $\bbR$. The scalar 
versions of $\rI^\rw_\vphi (\tha )$ and $\rI (\tha )$:
\beq\label{eq:wFIS}\beacl I^\rw_\vphi (\tha )&=\diy\int\vphi\,{\mathbf 1}(p_\tha >0)
p_{\tha}^{-1}\big((p_\tha)'\big)^2\\
\;&=\bbE_\tha\Big[\vphi (X){\mathbf 1}(p_\tha (X)>0)
p_\tha (X)^{-2}(p_\tha )'(X)^2 \Big],\\
I (\tha )&=\diy\int{\mathbf 1}(p_\tha >0)p_{\tha}^{-1}\big((p_\tha)'\big)^2\\
\;&=\bbE_\tha\Big[{\mathbf 1}(p_\tha (X)>0)
p_\tha (X)^{-2}(p_\tha )'(X)^2 \Big].\ena\eeq

For future use, we also introduce vectors:
\beq\label{eq:wFiV}\beacl V^\rw_\vphi (\tha )&=\diy\int\vphi\nabla p_\tha \\
\;&=\bbE_\tha\Big[\vphi (X){\mathbf 1}(p_\tha (X)>0)p_\tha (X)^{-1}
\nabla p_\tha (X)\Big],\ena\eeq
or -- when $d=1$ -- scalar values
\beq\label{eq:wFiVs}\beacl V^\rw_\vphi (\tha )&=\diy\int\vphi (p_\tha )'\\
\;&=\bbE_\tha\Big[\vphi (X){\mathbf 1}(p_\tha (X)>0)p_\tha (X)^{-1}
(p_\tha )'(X)\Big].\ena\eeq

Finally, set:
\beq E(\tha)=\diy\int\vphi p_\tha =\bbE_{\tha}[\varphi(X)]\eeq
and note that when $\vphi\equiv 1$ (a non-weighted case), $E(\tha )\equiv 1$.
$\Blt$\vskip .5cm

The above definitions lead to a useful extension of protagonist properties that are
well-known for the corresponding non-weighted quantities. As an example, consider a connection
between the weighted Kullback--Leibler divergence 
$${\tt K}^\rw_\vphi (p_\tha\| p_{\tha'})=\diy\int\vphi{\bf 1}(p_\tha >0)p_\tha\ln(p_\tha/ p_{\tha'})$$ 
(see \eqref{eq:wKuLe}) and the weighted Fisher information matrix $\rI^\rw_\vphi (\tha )$ from 
\eqref{eq:wFiM}. Fix a parameter value $\tha\in\Tha$; then, as $\tha'\to\tha$, one may expect that 
\beq\label{eq:relat}\beac\dfrac{1}{\tha_l -\tha'_l}\,{\tt K}^\rw_\vphi (p_\tha\| p_{\tha'})\to -E'_{\tha_l}(\tha ),\\
\dfrac{1}{(\tha_l -\tha'_l)(\tha_m-\tha'_m)}\Big[{\tt K}^\rw_\vphi (p_\tha\| p_{\tha'})\qquad\qquad\qquad{}\\
\qquad\qquad\qquad +\big(E(\tha )-E(\tha')\big)\Big]\to\dfrac{1}{2}I^\rw_{\vphi ,lm}(\tha).\ena\eeq
Here $E'_{\tha_l}(\tha )=\dfrac{\partial}{\partial\theta_l}E(\tha )$ and 
$I^\rw_{\vphi ,lm}(\tha)=\diy\int\vphi\,{\mathbf 1}(p_\tha >0)
p_{\tha}^{-1}(p_\tha)'_{\tha_l}(p_\tha)'_{\tha_m}$ stands for the entry $l,m$ of $\rI^\rw_\vphi$.
\vskip .5cm

To guarantee Eqns \eqref{eq:relat}, we need to assume that higher-order contributions vanish. 
A `physical' meaning of such an assumption is that for $\tha'$ close to $\tha$, the divergence
%${\tt K}_{\varphi}^\rw(p_{\theta}||p_{\theta'})$ 
admits a representation 
$$\beal 
{\tt K}_{\varphi}^\rw(p_{\theta}||p_{\theta'})=-\diy\int\varphi p_{\theta}\ln\dfrac{p_{\theta'}}{p_{\theta}}
=-\diy\int\varphi\,p_\tha \ln\Big[1+p_{\theta}^{-1}\nabla p_{\theta}(\theta'-\theta)^\rT\\
\qquad\qquad\qquad\;+\dfrac{1}{2}p_{\theta}^{-1}(\theta'-\theta)\big(\nabla p_{\theta}^T\nabla p_{\theta}\big)
(\theta'-\theta)^T+o(||\theta'-\theta||^2)\Big]\\
%=-\nabla E(\theta)(\theta'-\theta)^\rT+\dfrac{1}{2}(\theta'-\theta)\rI_{\theta}^\rw(\theta)(\theta'-\theta)^T
%+o(||\theta'-\theta||^2)
\ena$$
in which the gradients $\nabla$ and the remainder $o(||\theta'-\theta||^2)$ can be pulled out of integral
$\diy\int=\int_\cX$\, without affecting asymptotics under consideration. Such an approach is universally 
accepted in the literature and can be followed
in the weighted case as well.   \vskip .5cm

For the rest of this section we adopt the scheme  of $n$ IID observations where a random 
sample $\bX_1^n=(X_1,\ldots ,X_n)$ has been generated, with a joint PDF $p^{(n)}_\tha (\bx_1^n)
=\prod\limits_{i=1}^np_\tha (x_i)$. 
Our goal is to analyze an estimator $\tha^*(\bx_1^n)=(\tha^*_1(\bx_1^n),\ldots ,\tha^*_d(\bx_1^n))$ 
of parameter $\tha$ under a WF $\vphi^{(n)}$ written as a product:
\beq\label{eq:vphind}\vphi^{(n)}(\bx_1^n)=\prod\limits_{i=1}^n\vphi (x_i),\;\;\bx_1^n=(x_1,\ldots ,x_n)
\in\cX^n.\eeq
Here the map $x\in\cX\mapsto \vphi (x)\geq 0$ provides a WF in a single observation. 
\vskip .5cm

{\bf Definition.} Given an estimator $\tha^*(\bx_1^n)$,  set
\beq\label{eq:W=}
W(\tha ,\tha^*)=\diy\int\vphi^{(n)}\tha^*(\,\cdot\,)p^{(n)}_\tha
=\bbE^{(n)}_\tha\Big[\vphi^{(n)} (\bX_1^n)\tha^*(\bX_1^n)\Big].\eeq
Here, and below $\bbE^{(n)}_\tha$ means expectation under PDF $p^{(n)}$.
The value $W(\tha ,\tha^*)$ gives the weighted mean of $\tha^*$.

We split the value $W(\tha ,\tha^*)$ as follows:
\beq\label{wbias}\beal W (\tha,\tha^* )
=E(\tha)^n\tha + \bb(\tha),\;\;\tha\in\Theta .\end{array}\eeq
Here vector $\bb(\tha )=(b_1(\tha ),\ldots ,b_d(\tha ))$ (a scalar when $\Tha\subset\bbR$) 
can be interpreted as a weighted bias. \quad $\Blt$\vskip .5cm   

A natural assumption is that for a given $\tha\in\Theta$, $p_\tha$ and $\tha^*$ are regular, 
in the sense that the gradient and integration/expectation can be interchanged
\beq\label{eq:reg1}\beal
\nabla E(\tha)=\diy\int\vphi\nabla p_\tha,\\
\nabla W (\tha ,\tha^*)=\diy\int\vphi^{(n)}\tha^*(\,\cdot\,)^{\rT}\nabla p^{(n)}_\tha \\
\qquad\quad =\sum\limits_{i=1}^n\bbE^{(n)}_\tha\left[\vphi^{(n)}(\bX_1^n){\mathbf 1}(p_\tha (X_i)>0)
\tha^*(\bX_1^n)^{\rT}\dfrac{\nabla p_\tha (X_i)}{p_\tha(X_i)}\right].\ena\eeq
(The top equation is an identity for vectors while the bottom one is for matrices.)
Formally, assumption \eqref{eq:reg1} means that all involved gradients/derivatives exist
and the equalities hold true. A similar meaning is assigned to condition \eqref{eq:reg2}
below. 

Our goal in this section is to establish several forms of a weighted Cram\'er--Rao inequality for 
the weighted mean of the quadratic deviation of an estimator $\tha^*$ from the 
value of parameter $\tha$. An alternative view could be to consider a deviation 
from $E(\tha )^n\theta$; this line of work will be done elsewhere.  

Our arguments will involve 
the weighted Fisher information matrix for $n$ IID observations 
\beq\label{eq:wFInn}\beacl \rI^\rw_{\vphi^{(n)}}(\tha )&=nE(\tha )^{n-1}\rI^\rw_\vphi (\tha)\\
\;&\quad +n(n-1)E(\tha )^{n-2}V^\rw_\vphi (\tha )^\rT V^\rw_\vphi (\tha ).\ena\eeq
as well as its standard/non-weigted prototype with $\vphi\equiv 1$:
\vskip .5cm

{\bf Theorem 3.1.} (A weighted Cram\'er-Rao inequality, version (A)) {\sl Suppose that,
for a given $\tha\in\Tha$, matrix $\rI^\rw_\vphi (\tha )$ has finite diagonal entries
$I^\rw_{\vphi ,ll}(\tha )=\bbE_\tha\Big[\vphi (X){\mathbf 1}(p_\tha (X)>0)p_\tha (X)^{-2}
\left((p_\tha)'_{\tha_l}) (X)\right)^2\Big]$, $1\leq l\leq d$ (see \eqref{eq:wFiM}), 
and condition \eqref{eq:reg1} is fulfilled. Then the weighted squared quadratic deviation 
of estimator $\tha^*$ for $n\geq 2$ obeys the bound
\beq\label{weightKR}\beal
\bbE_\tha\Big[\vphi^{(n)}(\bX_1^n)\|\tha^*(\bX_1^n)-\tha\|^2\Big]\geq R(\tha ,n),\;\hbox{ where }\\
R(\tha,n)=\onwl{\max}\limits_{1\leq l\leq d}\Big[
\sum\limits_{k=1}^d\Big(E(\tha )^n+ (b_k)'_{\tha_l}(\tha )\Big)^2\\
\quad\times\big(nI^\rw_{\vphi ,ll}(\tha )E(\tha )^{n-1}
+n(n-1)(E'_{\tha_l}(\tha ))^2E(\tha )^{n-2}\big)^{-1}\Big],
\quad\tha\in\Tha .\ena\eeq

The equality is attained iff \ {\rPDF}s \ $p_\tha$, $\tha\in\Theta$, form
an exponential family: $p_\tha (x)=\exp\Big\{\nabla F(\tha )[t(x)^{\rT}-\tha^{\rT}]+F(\tha )+k(x)\Big\}$,
$x\in\cX$, and $\tha^*(\bX_1^n)=\dfrac{1}{n}\sum\limits_{i=1}^nt(X_i)$ is the sufficient statistic of the 
family.  
 
For $d=1$, bound \eqref{weightKR} becomes: 
\beq\label{wightKRs}\beal
\bbE_\tha\Big[\vphi^{(n)}(\bX_1^n)|\tha^*(\bX_1^n)-\tha)|^2\Big]
\geq R(\tha ,n)\;\;\hbox{where}\\
 R(\tha ,n)=\Big(E(\tha )^n+b'_{\tha}(\tha )\Big)^2\\
\qquad\times
\Big(nI^{\rw}_{\vphi} (\tha )E(\tha )^{n-1}+n(n-1)(E'(\tha ))^2E(\tha )^{n-2}\Big)^{-1},
\quad\tha\in\Tha ,\ena\eeq
and $I^\rw_\vphi (\tha )$ is the scalar version of matrix $\rI^\rw_\vphi (\tha )$; (see \eqref{eq:wFIS}).}
\vskip .5cm

{\bf Remark.} The RHS in \eqref{weightKR} takes the form 
$E(\tha )^nn^{-1}\onwl{\max}\limits_{1\leq l\leq d}\Big[I^\rw_{\vphi ,ll}(\tha )^{-1}\Big]$
when vector $\nabla E(\tha )=0$ and matrix $\nabla{\bf b}(\tha )=0$. Similarly, the RHS in 
\eqref{wightKRs} becomes $E(\tha )^nn^{-1}I^\rw_\vphi (\tha )^{-1}$ when $E'(\tha )=0$ and 
$b' (\tha )=0$. Both sets of conditions mean that point $\tha$ is a local extremum for $E(\tha )$
and ${\bf b}(\tha )$ or $b(\tha )$.
In this situation the weighted Cram\'er--Rao inequality yields a 
lower bound for the whole class of estimators $\tha^*$ satisfying the bottom condition in
\eqref{eq:reg1}. In general, one can think that the lower bound holds for a class of estimators 
where the derivatives of weighted bias ${\bf b}(\tha )$ or $b(\tha )$ are appropriately bounded.  

\vskip .5cm

{\bf Definition.} For an alternative version of the weighted Cram\'er-Rao inequality, 
we employ the functional  
\beq\label{eq:Ze}Z(\tha ,\tha^* )=\int\left(\vphi^{(n)}\right)^{1/2}p^{(n)}\tha^*(\,\cdot\,)^{\rT} 
=\bbE^{(n)}_\tha \left[\left(\vphi^{(n)}(\bX_1^n)\right)^{1/2}\tha^*(\bX_1^n)\right]\eeq 
and write
\beq\label{wbias2}\beal Z (\tha ,\tha^* )
=s(\tha)^n\tha^{\rT} + \bc(\tha)^{\rT},\;\hbox{where }
s(\tha )=\bbE_\tha\Big[\big(\vphi (X)\big)^{1/2}\Big],\quad\tha\in\Theta .\ena\eeq
Here $\bc(\tha )$ is a vector $(c_1(\tha ),\ldots ,c_d(\tha ))$. Cf. \eqref{eq:W=}, \eqref{wbias}.
\quad $\Blt$\vskip .5cm
 
Assume regularity conditions as follows. (Cf. \eqref{eq:reg1}.)
\beq\label{eq:reg2}\beal
0=\diy\int_\cX\nabla p_\tha (x),\\
\nabla Z (\tha ,\tha^*)=\diy\int_{\cX^n}\big(\vphi^{(n)} (\bx_1^n)\big)^{1/2}\tha^*(\bx_1^n)^{\rT}
\nabla p^{(n)}_\tha (\bx_1^n)\\
\qquad\quad=\sum\limits_{i=1}^n\bbE_\tha\left[\big(\vphi^{(n)} (\bX_1^n)\big)^{1/2}
{\mathbf 1}(p_\tha (X_i)>0)\tha^*(\bX_1^n)^{\rT}
\dfrac{\nabla p_\tha (X_i)}{p_\tha (X_i)}\right].\ena\eeq

\vskip .5cm

{\bf Theorem 3.2.} 
(A weighted Cram\'er-Rao inequality, version (B)) {\sl Suppose that, for a given $\tha\in\Tha$,
matrix $\rI (\tha )$ has finite diagonal entries $I_{ll}(\tha )=\bbE_\tha\Big[{\mathbf 1}(p_\tha (X)>0)
p_\tha (X)^{-2}\left((p_\tha)'_{\tha_l} (X)\right)^2\Big]$, $1\leq l\leq d$ (see \eqref{eq:sFIM}), 
and condition \eqref{eq:reg2} is fulfilled. Then the following bound holds true  for $n\geq 1$:
\beq\label{weightKRn}\beal
\bbE_\tha\Big[\vphi^{(n)} (\bX_1^n)\|\tha^*(\bX_1^n)-\tha\|^2\Big]\geq S(\tha ,n)\;\;\hbox{where}\\
S(\tha ,n)=\onwl{\max}\limits_{1\leq l\leq d}
\left[\sum\limits_{k=1}^d\Big(s(\tha )^n+(c_k)'_{\tha_l}(\tha )\Big)^2\big(nI_{ll}(\tha )\big)^{-1}\right].
\ena\eeq

For $d=1$ bound \eqref{weightKRn} becomes
\beq\beal\label{CR2}\bbE_\tha\Big[\vphi^{(n)}(\bX)|\tha^*(\bX)-\tha |^2\Big]
\geq\Big(s(\tha )^n+ c'(\tha )\Big)^2\big(nI (\tha )\big)^{-1}\ena\eeq
where $I(\tha )$ is a scalar version of $\rI (\tha)$ (see \eqref{eq:wFIS}).}

\vskip .5cm

{\it Proof of Theorems $3.1$ and $3.2$.} (A) Let us omit indices $1$ and $n$ in $\bX_1^n$. 
Combining 
\eqref{eq:reg1} and \eqref{wbias}, we get the matrix equality 
\beq\label{eq:gradE1}\bear
\sum\limits_{i=1}^n\bbE_\tha\Big[\vphi^{(n)}(\bX){\mathbf 1}(p_\tha (X_i)>0)
p_\tha (X_i)^{-1}\tha^*(\bX)^{\rT}\nabla p_\tha (X_i)\Big]\qquad{}\\
=nE(\tha )^{n-1}\tha^{\rT}\nabla E(\tha )
+E(\tha )^n \rU+\nabla \bb(\tha ).\ena\eeq
Here and below, $\rU$ stands for the unit $d\times d$ matrix. Next, the term 
$nE(\tha )^{n-1}\tha^{\rT}\nabla E(\tha )$ equals
$$\label{eq:gradE1A}\beal 
nE(\tha )^{n-1}\bbE_\tha\Big[\vphi (X_1){\mathbf 1}(p_\tha (X_1)>0)p_\tha (X_1)^{-1}\tha^{\rT}
\nabla p_\tha (X_1)\Big]\\
\qquad =\sum\limits_{i=1}^n\bbE_\tha\Big[\vphi^{(n)} (\bX){\mathbf 1}(p_\tha (X_i)>0)
p_\tha (X_i)^{-1}\tha^{\rT}\nabla p_\tha (X_i)\Big].\ena$$

Then \eqref{eq:gradE1} takes the form
\beq\label{eq:gradE2}\beal \sum\limits_{i=1}^n\bbE_\tha\bigg[\vphi^{(n)} (\bX)
{\mathbf 1}(p_\tha (X_i)>0)p_\tha (X_i)^{-1}(\tha^*(\bX)-\tha)^\rT\nabla p_\tha (X_i)\bigg]\\
=E(\tha )^n\rU+\nabla \bb(\tha ).\ena\eeq
%where $E$ stands for the $d\times d$ identity matrix.
Consequently, for an entry \ $k,l$ \ with $1\leq k,l\leq d$ \ we obtain
\beq\label{eq:gradE3}\beal
\bbE_\tha\bigg[\vphi^{(n)} (\bX)\Big(\tha^*_k(\bX)-\tha_k\Big)\sum\limits_{i=1}^n
{\mathbf 1}(p_\tha (X_i)>0)\\
\qquad\times p_\tha (X_i)^{-1}(p_\tha )'_{\tha_l}(X_i)\bigg]
=E(\tha )^n\delta_{k,l}+(b_k)'_{\tha_l}(\tha ).\ena\eeq

By the Cauchy--Schwarz inequality, 
\beq\label{CS1}\beal \left(E(\tha )^n\delta_{k,l}+(b_k)'_{\tha_l}(\tha )\right)^2
\leq\bbE_\tha\Big[\vphi^{(n)} (\bX)\big(\tha^*_k(\bX)-\tha_k)^2\Big]\\
\qquad\times\bbE_\tha\left[\vphi^{(n)} (\bX)\left(\sum\limits_{i=1}^n{\mathbf 1}(p(X_i)>0)
p_\tha (X_i)^{-1}(p_\tha)'_{\tha_l} (X_i)\right)^2\right].\ena\eeq
The factor \ $\bbE_\tha\left[\vphi^{(n)} (\bX)\left(\sum\limits_{i=1}^n{\mathbf 1}(p(X_i)>0)
p_\tha (X_i)^{-1}(p_\tha)'_{\tha_l} (X_i)\right)^2\right]$ \ equals
\beq\beal n\bbE_\tha\left[\vphi^{(n)}(\bX){\mathbf 1}(p(X_1)>0)p_\tha (X_1)^{-2}
\left((p_\tha)'_{\tha_l} (X_i)\right)^2\right]\\
\qquad + n(n-1)\left(\bbE_\tha\left[\vphi^{(n)}(\bX){\mathbf 1}(p(X_1)>0)
p_\tha (X_1)^{-1}(p_\tha)'_{\tha_l} (X_1)\right]\right)^2\\
\quad = nI_l(\tha )E(\tha )^{n-1}+n(n-1)\big(E'_{\tha_l}(\tha )\big)^2E(\tha)^{n-2}.\ena\eeq
Here we write $\bbE_\tha\left[\vphi (X_1){\mathbf 1}(p(X_1)>0)p_\tha (X_1)^{-1}(p_\tha)'_{\tha_l} (X_1)
\right]=E'_{\tha_l}(\tha )$ 
by virtue of the top equation in \eqref{eq:reg1}.

Therefore, 
$$\beal\bbE_\tha\Big[\vphi^{(n)} (\bX)\big(\tha^*_k(\bX)-\tha_k)^2\Big]
\geq\Big(E(\tha )^n\delta_{k,l}+(b_k)'_{\tha_l}(\tha )\Big)^2\\
\qquad\qquad\times\Big(nI_l(\tha )E(\tha )^{n-1}+n(n-1)(E'_{\tha_l}(\tha ))^2E(\tha )^{n-2}\Big)^{-1},\ena$$
and summing over $k$ yields
$$\beal
\bbE_\tha\Big[\vphi^{(n)} (\bX)\|\tha^*(\bX)-\tha\|^2\Big]
\geq\sum\limits_{k=1}^d\Big(E(\tha )^n+(b_k)'_{\tha_l}(\tha )\Big)^2\\
\qquad\qquad\times\Big(nI_l(\tha )E(\tha )^{n-1}+n(n-1)(E'_{\tha_l}(\tha ))^2E(\tha )^{n-2}\Big)^{-1}.\ena$$
Optimizing over \ $l$ \ leads to \eqref{weightKR}. For equality in the bound we 
need equality in the Cauchy-Schwarz inequality, which leads to an exponential family.\vskip .5cm

(B) We can repeat the argument where $\vphi^{(n)}$ is replaced by $(\vphi^{(n)})^{1/2}$.
The starting point is the equality similar to \eqref{eq:gradE1}: 
\beq\label{eq:gradE11}\bear
\sum\limits_{i=1}^n
\bbE_\tha\left[\big(\vphi^{(n)}  (\bX)\big)^{1/2}\tha^*(\bX)^{\rT}{\mathbf 1}(p(X_i)>0)
p_\tha (X_i)^{-1}\nabla p_\tha (X_i)\right]\qquad{}\\
\qquad\qquad\qquad\qquad =ns(\tha )^{n-1}\tha^{\rT}\nabla s(\tha )
+s(\tha )^n\rU+\nabla \bc(\tha ).\ena\eeq
This is followed by
$$\beal ns(\tha )^{n-1}\tha^{\rT}\nabla s(\tha )\\
\quad =ns(\tha )^{n-1}\bbE_\tha\left[\big(\vphi (X_1)\big)^{1/2}\tha^{\rT}
{\mathbf 1}(p(X_1)>0)p_\tha (X_1)^{-1}\nabla p_\tha (X_1)\right]\\
\quad =\bbE_\tha\left[\big(\vphi^{(n)}  (\bX)\big)^{1/2}\tha^{\rT}\sum\limits_{i=1}^n
{\mathbf 1}(p(X_i)>0)p_\tha (X_i)^{-1}\nabla p_\tha (X_i)\right].\ena$$
Then \eqref{eq:gradE11} takes the form
$$\label{eq:gradE22}\bear
\bbE_\tha\left[\big(\vphi^{(n)} (\bX)\big)^{1/2}\Big(\tha^*(\bX)-\tha\Big)^{\rT}
\sum\limits_{i=1}^n{\mathbf 1}(p(X_i)>0)
\dfrac{\nabla p_\tha (X_i)}{p_\tha (X_i)}\right]\quad{}\\
=s(\tha )^n\rU+\nabla \bc(\tha ).\ena$$

Again, it is  convenient to write down the equation for a single entry:
$$\label{eq:gradE33}\bear
\bbE_\tha\left[\big(\vphi^{(n)} (\bX)\big)^{1/2}\Big(\tha^*_k(\bX)-\tha_k\Big)\sum\limits_{i=1}^n
{\mathbf 1}(p(X_i)>0)p_\tha (X_i)^{-1}(p_\tha)'_{\tha_l} (X_i)\right]\qquad{}\\
=s(\tha )^n\delta_{k,l}+(c_k)'_{\tha_l}(\tha ).\ena$$
By Cauchy--Schwarz, 
$$\label{CS11}\beal \left(s(\tha )^n\delta_{k,l}+(c_k)'_{\tha_l}(\tha )\right)^2
\leq\bbE_\tha\Big[\vphi^{(n)} (\bX)\big(\tha^*_k(\bX)-\tha_k)^2\Big]\\
\qquad\qquad\times\bbE_\tha\left[\left(\sum\limits_{i=1}^n{\mathbf 1}(p(X_i)>0)
p_\tha (X_i)^{-1}(p_\tha)'_{\tha_l} (X_i)\right)^2\right].\ena$$

Next, the term $\bbE_\tha\left[\left(\sum\limits_{i=1}^n{\mathbf 1}(p(X_i)>0)
p_\tha (X_i)^{-1}(p_\tha)'_{\tha_l} (X_i)\right)^2\right]$ equals
\beq\beal n\bbE_\tha\left[\left({\mathbf 1}(p(X_1)>0)p_\tha (X_1)^{-1}(p_\tha)'_{\tha_l} (X_1)\right)^2\right]\\
\quad + n(n-1)\left(\bbE_\tha\left[{\mathbf 1}(p(X_1)>0)p_\tha (X_1)^{-1}(p_\tha)'_{\tha_l} (X_1)
\right]\right)^2 = nI_l(\tha ).\ena\eeq
Here we write $\bbE_\tha\left[{\bf 1}(p_\tha (X_1)>0)p_\tha (X_1)^{-1}(p_\tha)'_{\tha_l} (X_1)\right]=0$,
owing to the top equation in \eqref{eq:reg2}.

Consequently, 
$$\beal\bbE_\tha\Big[\vphi^{(n)} (\bX)\big(\tha^*_k(\bX)-\tha_k)^2\Big]
\geq\Big(s(\tha )^n\delta_{k,l}+(c_k)'_{\tha_l}(\tha )\Big)^2(nI_l(\tha ))^{-1},\ena$$
and summing over $k$ yields
$$\beal
\bbE_\tha\Big[\vphi^{(n)} (\bX)\|\tha^*(\bX)-\tha\|^2\Big]
\geq\sum\limits_{k=1}^d\Big(s(\tha )^n+(c_k)'_{\tha_l}(\tha )\Big)^2(nI_l(\tha ))^{-1}.
\ena$$
Now, optimizing over \ $l$ \ leads to \eqref{weightKRn}.  $\square$ \vskip .5cm
\vskip .5cm

Examples of bounds emerging from Theorems 3.1 and 3.2 are discussed in Appendix C.\vskip .5cm

{\bf Remarks. 1.} Equality in Theorem 3.2 in general is not attainable. A more detailed analysis
will be given in a separate work.\vskip .5cm

{\bf 2.} The Cauchy--Schwarz inequality can be also used for each term in \eqref{eq:gradE3}:
$$\beal E(\tha )^n\delta_{k,l}+(b_k)'_{\tha_l}(\tha )\leq 
\left(\bbE_\tha\Big[\vphi^{(n)} (\bX)\big(\tha^*_k(\bX)-\tha_k)^2\Big]\right)^{1/2}\\
\quad\times\sum\limits_{i=1}^n
\left(\bbE_\tha\left[\vphi^{(n)} (\bX) {\mathbf 1}(p_\tha (X_i)>0)
\left(p_\tha (X_i)^{-1}(p_\tha)'_{\tha_l} (X_i)\right)^2\right]\right)^{1/2}\ena$$
or
$$\beal\bbE_\tha\Big[\vphi^{(n)} (\bX)\big(\tha^*_k(\bX)-\tha_k)^2\Big]\\
\qquad\qquad\geq\Big(E(\tha )^n\delta_{k,l}+(b_k)'_{\tha_l}(\tha )\Big)^2 
\left(n I_l(\tha )^{1/2}\right)^{-2}.\ena$$
However, it leads to a bound with a factor $n^{-2}$, and we will not discuss it in this paper.
\vskip .5cm

{\bf 3.} Assume that domain $\Tha =\bbR^d$. Let $\tha^*$ be an unbiased 
estimator of parameter $\tha$, with vector
${\bf b}(\tha )=0$ $\forall$ $\tha\in\Tha$; cf.\eqref{wbias}.
Consider the weighted covariance between $\tha^*$ and the random vector 
$\ups =\ups (\bX_1^n):=\nabla\ln p^{(n)}_\tha (\bX_1^n)$:
\beq\label{eq:CovR}{\rm Cov}^\rw_{\vphi}(\ups, \tha^*):=\bbE_{\tha}[\vphi^{(n)}(\bX_1^n)
\ups (\bX_1^n)\tha^*(\bX_1^n)].\eeq
Then the following two facts (i), (ii) are equivalent. (i) The covariance \eqref{eq:CovR}
depends on $\tha$ through $E(\tha )$ only. (ii) For any estimator $\widetilde\tha$ for which 
$W(\tha ,{\widetilde\tha})=0$ (see \eqref{eq:W=}), the weighted covariance 
${\rm Cov}^\rw_{\vphi}(\ups ,{\widetilde\tha})=0$. 
\vskip.5cm

To prove that (i) implies (ii), note that 
the weighted mean of estimators $\tha^*$ and $\tha^*+{\tilde\tha}$ are equal. However,
$${\rm Cov}^\rw_{\vphi}(\ups ,\tha^*+{\tilde\tha})={\rm Cov}^\rw_{\vphi}(\ups ,\tha^*)
+{\rm Cov}^\rw_{\vphi}(\ups ,{\widetilde\tha})={\rm Cov}^\rw_{\vphi}(\ups ,\tha^*)$$ 
which implies that ${\rm Cov}^\rw_{\vphi}(\ups ,{\widetilde\tha})=0$.

Conversely, if the weighted covariance between $\ups$ and ${\widetilde\tha}$ vanishes  
then for any two estimators with
$\bbE_{\tha}[\vphi^{(n)}(\bX_1^n)\tha_1^*(\bX_1^n)]
=\bbE_{\tha}[\vphi^{(n)}(\bX_1^n)\tha_2^*(\bX_1^n)]$ we have 
${\rm Cov}_{\vphi}(\ups ,\tha_1^*-\tha_2^*)=0$ as $\tha_1^*-\tha_2^*$ 
has $W(\tha ,\tha^*_1-\tha^*_2)=0$. 
\vskip .5cm

We now pass to weighted van Trees inequalities.
For motivation for and applications of such inequalities in a non-weighted case, see \cite{vT};
proofs and extensions are discussed in \cite{GiLe}. In Theorems 3.3 and 3.4 below we assume 
that parameter 
$\tha$ has a prior distribution $\Pi (\rd\tha )$ on $\Theta$ and refer to the quantities $R(\tha ,n)$ and 
$S(\tha ,n)$ are determined in \eqref{weightKR} and \eqref{weightKRn}.  
\vskip  .5cm

{\bf Theorem 3.3.} (Weighted van Trees inequality, version A) {\sl Under condition 
\eqref{eq:reg1}, for any probability distribution $\Pi$ on $\Theta$  
\beq\label{weightvT1}\beal
\diy\int_\Theta\bbE_\tha\Big[\vphi^{(n)} (\bX_1^n)\|\tha^*(\bX_1^n)-\tha\|^2\Big]\Pi (\rd\tha)
\geq\diy\int_\Theta R(\tha ,n)\Pi (\rd\tha) .\ena\eeq}

{\bf Theorem 3.4.} (Weighted van Trees inequality, version B) {\sl Under condition 
\eqref{eq:reg2}, for any probability distribution $\Pi$ on $\Theta$  
\beq\label{weightvT2}\beal
\diy\int_\Theta\bbE_\tha\Big[\vphi^{(n)} (\bX_1^n)\|\tha^*(\bX_1^n)-\tha\|^2\Big]\Pi (\rd\tha)
\geq\diy\int_\Theta S(\tha ,n)\Pi (\rd\tha).\ena\eeq}

\vskip .5cm

{\it Proof.} \quad Inequalities \eqref{weightvT1} and \eqref{weightvT2} are direct corollaries of 
Proposition 2.14. $\square$
\vskip .5cm 

In Theorem 3.4  we use condition \eqref{eq:reg3}, \eqref{eq:Tea}, involving a smooth PDF 
$\pi$ on $\Theta =\bbR^d$. (The case where $\Theta$ is a torus or another closed manifold 
can be covered as well, at the expense of additional technicalities.) 
Consider the weighted Fisher information density for parameter $\tha$:
$$j^\rw_\varphi (\tha ):={\bf 1}(\pi (\tha )>0)E(\tha)^n
\dfrac{\|\nabla\pi (\tha)\|^2}{\pi(\tha)},\quad\tha\in\bbR^d.$$
Recall, $I (\tha)$ stands for Fisher information matrix in a given family  
$p_\tha =p(\,\cdot\,,\tha )$; cf. \eqref{eq:sFIM}.
 
Consider the following condition: \ $\forall$ \ $n=1,2,\ldots$,
\beq\label{eq:reg3}\beal
\lim\limits_{h\to{\bf 0}}\dfrac {1}{\|h\|^2}\diy\int_{\bbR^d}{\bf 1}(\pi (\tha )>0)\left(\diy\int\vphi^{(n)} 
{\bf 1}\big(p^{(n)}_\tha >0\big)\right.\\
\qquad\qquad\times\left.\dfrac{\big[\pi(\tha+h)p^{(n)}_{\tha+h}-\pi(\tha)p^{(n)}_\tha\big]^2
}{\pi(\tha)p^{(n)}_\tha}\right){\rd}\tha =T(\vphi ,\pi)\ena\eeq
%\diy\int_\Theta E(\tha)^n{\mathbf 1}(\pi (\tha >0))\dfrac{||\nabla_{\tha}\pi(\tha||^2}{\pi(\tha)}%\rd\tha+
where
\beq\label{eq:Tea} T(\vphi ,\pi)=\diy\int_{\bbR^d}\Big[j^\rw_\vphi (\tha )+\big(n\,{\rm tr}\,[ \rI^\rw_\varphi (\tha)]
+n(n-1)\|\nabla E(\tha )\|^2\big)\pi (\tha )\Big]{\rd}\tha .\eeq

The meaning of condition \eqref{eq:reg3} is that the limit $\lim\limits_{h\to{\bf 0}}$ can be
performed under the integrals. After the change of the order of operations one comes to the
integral 
$$\beal\diy\int_{\bbR^d}{\bf 1}(\pi (\tha )>0)\diy\left(\int\vphi^{(n)}{\mathbf 1}(p(\,\cdot\,,\tha)>0)
\dfrac{\|p^{(n)}_\tha\nabla\pi (\tha )+\pi (\tha )\nabla p^{(n)}_\tha \|^2}{
\pi (\tha )p^{(n)}_\tha}\right)\rd\tha\ena$$
which leads to the value $T(\vphi ,\tha )$ in \eqref{eq:Tea}. Admittedly, a verification of 
condition \eqref{eq:reg3} (as well as \eqref{eq:reg2} and \eqref{eq:reg1}) can be of an
academic interest only, as in applications it is routinely assumed to hold.
\vskip .5 cm

{\bf Theorem 3.5.} (Weighted van Trees inequality, version C) {\sl Under condition \eqref{eq:reg3},}
\beq\label{weightvT3}\bear
\diy\int_{\bbR^d} \bbE_{\tha}\Big[\vphi^{(n)}(\bX_1^n)\big\|\tha^*(\bX_1^n)-\tha\big\|^2\Big]
\pi(\tha){\rd}\tha\qquad\qquad\qquad{}\\ 
\geq\Big[\diy\int_{\bbR^d} E(\tha)^n \pi(\tha){\rd}\tha\Big]^2\,(T(\vphi ,\pi))^{-1}.
\ena\eeq
{\sl where $T(\vphi ,\tha )$ is given by the \rRHS \ in \eqref {eq:Tea}.}
\vskip .5cm

{\it Proof.} As before, write $\bx$ and $\bX$ in place of $\bx_1^n$ and $\bX_1^n$.
Consider a family $\{\gam_h,\,h\in\bbR^d\}$ of joint PDFs on $\cX^n\times\bbR^d$:
\beq\label{eq:gamma}
\gamma_h(\bx,\tha)=\pi(\tha+h)p(\bx,\tha+h),\quad \bx\in {\cal X}^n, \tha\in \bbR^d.
\eeq
To avoid confusion write ${\cal E}_h$ for integrals in both variables $\bx$, $\tha$. Set:
\beq
D_h(\bx,\tha)=\dfrac{\gamma_h(\bx,\tha)-\gamma_0(\bx,\tha)}{\gamma_0(\bx,\tha)}
\eeq
and 
\beq\beal
G_{\varphi}(h):={\cal E}_h[\vphi^{(n)}(\bX)(\tha^{\rT}(\bX)-\tha^{\rT})],
\quad\hbox{with}\\
G_{\varphi}(h)=G_{\varphi}(0)+ h^{\rT}\diy\int E(\tha +h)^n\pi(\tha+h){\rd}\tha\\
\qquad\quad=G_{\varphi}(h)=G_{\varphi}(0)+ h^{\rT}\diy\int E(\tha)^n\pi(\tha){\rd}\tha ,\quad h\in\bbR^d.
\ena\eeq
Clearly,
\beq\beal
{\cal E}_0[\vphi^{(n)}(\bX)D_h(\bX,\tha)(\tha^*(\bX)^{\rT}-\tha^{\rT})]\\
\qquad={\cal E}_h[\vphi^{(n)}(\bX)(\tha^*(\bX)^{\rT}-\tha^{\rT})]
-{\cal E}_0[\vphi^{(n)}(\bX)(\tha^*(\bX)^{\rT}-\tha^{\rT})]\\
\qquad=G_{\varphi}(h)-G_{\varphi}(0)=h^{\rT}\diy\int E(\tha)^n\pi(\tha){\rd}\tha .
\ena\eeq

Now use the Cauchy-Schwarz inequality: 
\beq\label{eq:CS12}\beal
||h||^2\left(\diy\int E(\tha)^n\pi(\tha){\rd}\tha\right)^2 
=\left\|{\cal E}_0\Big[\varphi(\bX)D_h(\bX,\tha)(\tha^*(\bX)^{\rT}-\tha^{\rT})
\Big]\right\|^2\\
\qquad\leq \left({\cal E}_0[\vphi^{(n)}(\bX)D_h(\bX, \tha)^2]\right)
\left({\cal E}_0\Big[\vphi^{(n)}(X)\|\tha^*(\bX)-\tha\|^2\Big]\right).\ena\eeq

The factor ${\cal E}_0\Big[\vphi^{(n)}(X)||\tha^*(\bX)-\tha||^2\Big]$ in the RHS of 
\eqref{eq:CS12} equals
\beq
\diy\int \pi(\tha)\bbE_{\tha}[\vphi^{(n)}(\bX)||\tha^*(\bX)-\tha||^2]{\rd}\tha
\eeq
giving the expression in the LHS of \eqref{weightvT3}. Next, owing to \eqref{eq:reg3},
\beq\beal
\dfrac{1}{\|h\|^2}{\cal E}_0[\vphi^{(n)}(\bX)D_h(\bX, \tha)^2]\to T(\vphi ,\pi)\quad\hbox{as $h\to{\bf 0}$,}
\ena\eeq
which implies \eqref{weightvT3}. $\square$
\vskip .5cm

{\bf Remark 3.} The assertion of Theorem 3.5 does not use any `regularity' assumption upon 
estimator $\tha^*$.  \vskip .5cm

\section{ Exponential families: formulas for divergences}

As we have seen in Lemma 2.2, the weighted loss is specified or bounded in terms of weighted
divergences. In the case of an exponential family, the related calculations are considerably simplified.
The machinery of these calculations involves some facts and notions which we discuss in this section.
\vskip .5cm

{\bf Definition.} Let $\Theta$ be a convex domain in $\bbR^d$. A family of PDFs 
$p_\tha$ on $\cX$, parametrized by $\tha =(\tha_1,\ldots ,\tha_d)\in \Theta$, is called exponential if 
\beq
p_\tha (x)=\exp\Big[\tha t(x)^\rT-F(\tha)+k(x)\Big],\quad x\in\cX.
\eeq
Function $t:\;\cX\to\bbR^d$ gives a sufficient statistic for the family. Function $F:\,\Theta\to
\bbR$, is called the log-normalizer.\footnote{In this section we use a restricted form of an exponential
PDF family in order to simplify the technical aspects of formulas below.} 

%Here $P_\tha$ stands for the probability distribution generated by PDF $p(\bx ,\tha )$.

Given a WF $\vphi $, define an adjoint exponential family $\{p^*_\tha ,\;\tha\in\Theta\}$, with
\beq
p_\tha^*(x)=\re^{\langle \tha,t^*(x)\rangle-F^*(\tha)+k^*(x)},\quad x\in\cX.\eeq
Here $t^*(x)=t(x)$,  $k^*(x)=k(x)+\ln\vphi (x)$, and $F^*(\tha)$ is the adjoint log-normalizer
\beq\label{lognor}
F^*(\tha)=\ln \diy\int\re^{\langle t(\,\cdot\,),\tha \rangle+k^*(\,\cdot\,)}
=F(\tha )+\ln\,E_\vphi (\tha ).\eeq
As before, $E_\vphi (\tha )=\diy\int\vphi\,p_\tha =\bbE_\tha [\vphi (X)]$ stands for the mean value of the weight
under PDF $p_\tha$.

Suppose that an exponential PDF family $\{p_\tha,\;\tha\in\Tha\}$ has been given. We recall that
the Bregman divergence (see \cite{Nie1}) is defined as 
\beq\label{eq:BreNW}
{\tt B}_F(\tha',\tha):=F(\tha')-F(\tha)-\langle\tha'-\tha, \nabla F(\tha)\rangle\eeq
and coincides with the standard Kullback--Leibler divergence:
\beq\label{eq:BKLNW}{\tt B}_F(\tha',\tha)={\tt K}(p_{\tha}||p_{\tha'}).\eeq
Next, the quantity
\beq\label{eq:eq:BreW}\beal
{\tt B}^\rw_{\vphi, F}(\tha',\tha):=\re^{F^*(\tha)-F(\tha)}\Big[F(\tha')-F(\tha)
-\langle\tha'-\tha, \nabla F^*(\tha)\rangle\Big]\\
\qquad\quad =E_\vphi (\tha)\Big[F(\tha')-F(\tha)-\langle (\tha'-\tha )^{\rT}
\nabla\Big(F(\tha)+\ln\,E_\vphi (\tha)\Big)\rangle\Big].\ena\eeq
is called the weighted Bregman divergence. \quad $\Blt$\vskip .5cm

{\bf Proposition 4.1.} {\sl The weighted Bregman divergence coincides with the weighted 
Kullback--Leibler divergence: 
\beq\label{eq:BKLW}{\tt B}^\rw_{\vphi, F}(\tha',\tha)={\tt K}^\rw_{\vphi}(p_{\tha}||p_{\tha'}).\eeq}
Cf. \eqref{eq:BKLNW}. \vskip .5cm

{\it Proof.}  By the direct differentiation
$$\nabla F^*(\tha)=\diy\int  t(\,\cdot\,) p_\tha^*
=\re^{F(\tha)-F^*(\tha)}\int t(\,\cdot\,)\vphi\, p_\tha .$$
It implies:
\beq\beacl{\tt K}^\rw_{\vphi}(p_{\tha}||p_{\tha'})&=\diy\int \vphi {\bf 1}(p_\tha >0)
p_\tha\ln\dfrac{p_\tha}{p_{\tha'}}  \\
\;&=\diy\int \vphi {\bf 1} (p_\tha >0) p_\tha \big[F(\tha')-F(\tha)-\langle\tha'-\tha, t(\,\cdot\,)\rangle\big]  \\
\;&=\re^{F^*(\tha)-F(\tha)}\Big[F(\tha')-F(\tha)-\langle\tha'-\tha, 
\nabla F^*(\tha)\rangle\Big]. \label{2.30}\end{array}\eeq
$\Box$
\vskip .5cm

Our next goal is to express weighted entropies and divergences for PDFs $p_\tha$ from an exponential 
family in terms of log-normalizers $F(\tha)$ and $F^*(\tha)$. For the weighted Shannon
entropy (see \eqref{eq2:wSE}) we have the expression 
\beq\label{2.344}\beacl
h^\rw_{\vphi}(p)&=E_\vphi (\tha ) \left[F(\tha)-\tha\cdot\nabla F^*(\tha)\right]\\
\;&=E_\vphi (\tha)[F(\tha)-\tha\cdot\nabla F(\tha)]-\tha\cdot\nabla E_\vphi (\tha).
\end{array} \eeq

Next, consider the quantity %weighted Jensen difference divergence 
\beq\label{eq:Jphi}\beacl
J_{\vphi ,\alpha}(p_\tha)&=\diy\int\vphi p^{\alpha}_\tha=\diy\int \re^{\langle t(\,\cdot\,), \alpha \tha\rangle -\alpha 
F(\tha)+\alpha k(\,\cdot\,)+ \ln(\vphi)}  \\
\;&= \re^{F(\alpha\tha)-\alpha F(\tha)}\diy\int \vphi p_{\alpha\tha}\re^{(\alpha-1)k(\,\cdot\,)}. \ena\eeq
It gives an explicit expression for the weighted R\'enyi entropy (cf. \eqref{RenyiEnt}):  
\beq\label{eq:REnt1}\beacl
r^\rw_{\vphi ,\alpha}(p_\tha)&=\dfrac{E_\vphi (\tha )}{1-\alpha}\ln\dfrac{J_{\vphi ,\alpha}(p_\tha)}{E_\vphi (\tha )}
=\dfrac{E_\vphi (\tha) }{1-\alpha}\bigg[F(\alpha\tha)-\alpha F(\tha)\\
\;&\quad+\ln\diy\int\vphi p_{\alpha\tha}
\re^{(\alpha -1)k(\,\cdot\,)}-F^*(\tha )+F(\tha )\bigg].\ena\eeq
When  %the so-called natural exponential models with 
$k(x)=0$  (e.g., in the case of an exponential, Gamma or Gaussian distribution), the formula simplifies:
\beq\beacl r^\rw_{\vphi ,\alpha}(p)&
=\dfrac{E_\vphi (\tha )}{1-\alpha}\Big[F^*(\alpha \tha)-\alpha F(\tha)-F^*(\tha )-F(\tha)\Big]\\
\;&=\dfrac{E_\vphi (\tha )}{1-\alpha}\Big[\ln E_\vphi (\alpha\tha)+F(\alpha\tha)-\alpha F(\tha)-\ln E_\vphi (\tha )\Big].
\ena\label{2.33}\eeq

As $\alpha\to 1$, applying the L'Hopital rule to the RHS of \eqref{2.33}, one gets
%for any WF $\vphi$, 
\beq\lim\limits_{\alpha\to 1}r^\rw_{\vphi, \alpha}(p)=-\diy\int_{\cX}\vphi (x)\Big[t(\bx)\tha^\rT
-F(\tha)\Big]p_\tha (\bx)  =h^\rw_{\vphi}(p).\label{2.49}\eeq
Indeed, the last term in the final line of ({\ref{2.33}}) does not depend on $\alpha$, and differentiating
$\ln\,E(\alpha\tha)$ gives $t(\bx)\tha^\rT-\nabla F(\alpha\tha)\tha^\rT$. Here one 
needs an assumption
that derivative in $\alpha$ can be carried under the integral. Finally, the term 
$-\nabla F(\alpha\tha)\tha^\rT$ cancels as a result of differentiation of $F(\alpha\tha )$
in the RHS of \eqref{2.33}. Then Eqn \eqref{2.49} follows. \vskip .5cm

Furthermore, for the weighted Chernoff coefficient (see \eqref{ChernCo}) we get  
\beq\beacl
\ln C^\rw_{\vphi ,\alpha}(p_\tha ,p_{\tha'})&=\ln \diy\int \vphi\,p_\tha^\alpha p_{\tha'}^{1-\alpha} 
-\ln\int\vphi p_\tha\\ 
\;&=F(\alpha\tha+(1-\alpha)\tha')-\alpha F(\tha)-(1-\alpha)F(\tha')\\
\;&\quad+\ln E_\vphi (\alpha\tha+(1-\alpha)\tha')-\ln\,E_\vphi (\tha ).
\end{array}\label{4.36}\eeq
Equivalently, one can use 
the weighted Burbea-Rao divergence ${\tt U}^\rw_{\alpha, F}(\tha,\tha')$, cf. \cite{BR}, \cite{NB}. 
Here
\beq
\beacl{\tt U}^\rw_{F,\alpha}(\tha,\tha')&:=\alpha F(\tha)+(1-\alpha)F(\tha')-F(\alpha\tha+(1-\alpha)
\tha')\\
\;&\,={\tt U}^\rw_{F,1-\alpha}(\tha',\tha ).\ena\eeq
Then the following identity emerges:
\beq\beal
C^\rw_{\vphi ,\alpha}(p_\tha ,p_{\tha'}) \\
\quad =\diy\int \vphi\re^{\langle t(\,\cdot\,),  \alpha\tha+(1-\alpha)\tha' \rangle
-(\alpha F(\tha)+(1-\alpha)F(\tha'))+k(\,\cdot\,)}\Big/\int\vphi p_\tha \\
\quad =\re^{F(\alpha\tha+(1-\alpha)\tha')-\alpha F(\tha)-(1-\alpha)F(\tha')} 
\dfrac{E_\vphi (\alpha\tha+(1-\alpha)\tha')}{E_\vphi (\tha )}\\
\quad =\re^{-{\tt U}^\rw_{F,\alpha}(\tha,\tha')}\dfrac{E_\vphi (\alpha\tha+(1-\alpha)\tha')}{E_\vphi (\tha )}.
\end{array}
\eeq 
Finally, for the Chernoff divergence (see \eqref{1.3.11}) we obtain
\beq\label{eq:4.18}{\tt C}^\rw_{\vphi ,\alpha}(p _{\tha},p_{\tha'})
={\tt U}^\rw_{F,\alpha}(\tha,\tha')-\ln E_\vphi\big(\alpha\tha+(1-\alpha)\tha'\big)
+\ln E_\vphi (\tha ).\eeq
\vskip .5cm

\section{Concluding remarks}

In the current paper we attempted to outline possible routes for the development of 
basic elements of context dependence in a statistical theory.
Formally, incorporating a weight function $\vphi$ should not change fundamentals of statistical inference;
all staple principles must have their weighted counterparts/extensions. We discussed a host of
properties that can be presented in an accessible form; a variety of facts and concepts 
collected in this paper makes first steps towards a new doctrine.

Let us outline examples of possible further steps.\footnote{In this paper we discuss basics of 
hypothesis testing and do not dwell on topics in other directions of Statistics, Probability or Information
theory.} (I) Consider the 
problem of constructing a uniformly 
most powerful weighted test
in an exponential family with a sufficient statistic $t(\,\cdot\, )$. If parameter $\tha$ has dimension $1$ 
then, regardless of the dimension  
of the data, such a test can be constructed, in terms of the mean value ${\mathbb E}_\tha[\vphi (X)t(X)]$
that is a monotone function of $\tha$.  It would be instructive to provide technical details and examples 
involving a suitable weighted power of a test. (II) A similar problem can be considered for families
with (again suitably defined) monotone likelihood ratio. (III) Weighted analogs of significance/goodness
of fit tests may lead to interesting generalizations of Pearson's $\chi^2$-, Student's $t$- and 
Fisher--Snedecor's $F$-distributions
and/or tests. (IV) Likewise, a weighted extension of the Kolmogorov--Smirnov test might emerge
and prove to be useful. In general, the context-dependent development of Non-parametric 
Statistics seems promising, particularly in connection with Information theory.  For a modern treatment of the 
standard Kolmogorov--Smirnov test and other aspects of non-parametric Statistical theory 
and its relations with Information theory and Coding, see \cite{GN}. %Sect 6.2.1.  
\vskip .5cm

{\bf Acknowledgment.} This work has been supported by the RSF (project No. 23-21-00052). 
%and the HSE University Basic Research Program. 
The support from PSU Math Department and St John's College,
Cambridge, is also acknowledged.

\newpage

\section{\bf Appendix A. Weighted TV distance for Gaussian PDFs}\vskip 1cm

In this section we specify the weighted TV distance between PDFs 
$p(x)=\dfrac{1}{\sqrt{2\pi}}\re^{-x^2/2}$ and 
$q(x)=\dfrac{1}{\sqrt{2\pi}}\re^{-(x-a)^2/2}$ where $a\geq 0$. The answers 
involve the function ${\rm Erf}(z)=\dfrac{2}{\sqrt{\pi}}\diy\int_0^z \re^{-t^2}{\rd}t$.\vskip .5cm

{\bf Example A1.}  Set $\cX=\mathbb R$ and $\vphi (x)=x^2+bx+c$ with $b\in\mathbb R$ and 
$c\geq b^2/4$. Then
$$
\tau^\rw_{\vphi}(p,q)=\dfrac{1}{2}\left[\sqrt{\dfrac{2}{\pi}}a 
\re^{-a^2/8}+a(a+b)+(2+a^2+ab+2c) {\rm Erf}\left(\dfrac{a}{2\sqrt{2}}\right)\right].$$
\vskip .5cm

{\bf Example A2.} Set: $\cX=\mathbb R$ and $\vphi (x)=|x|$. Then
$$
\tau^\rw_{\vphi}(p,q)=\dfrac{a}{2}\left(1+ {\rm Erf}\left(\dfrac{a}{2\sqrt{2}}\right)\right).$$
\vskip .5cm

{\bf Example A3.} Set: $\cX=\mathbb R$ and $\vphi (x)=\re^{\gam x}$, $\gam\in\mathbb R$.
Then
$$\beac\tau^\rw_\vphi (p,q)=
\re^{(\gam^2+2\gam a)/2}\left[{\rm{Erf}}\left(\dfrac{a+2\gam}{2\sqrt 2}\right)
+\re^{-\gam a}{\rm{Erf}}\left(\dfrac{a-2\gam}{2\sqrt 2}\right)-2\right].\ena$$

%=\re^{\lam^2/2}\Big[2\Phi\left(\lam+\dfrac{a}{2}\right)-1\Big]
%+\re^{-\lam a+\lam^2/2}\Big[2\Phi\left(-\lam +\dfrac{a}{2}\right)-1\Big]\ena\eeq
\vskip 1cm

\section {Appendix B. Examples of weighted Cram\'er--Rao bounds} 

{\bf{(B1)}} (A shift-parametric normal family) Consider a Gaussian PDF family 
$p_\tha =p_\mu$, $\mu\in\bbR$, where 
$$p_\mu (x)=\dfrac{1}{\sqrt{2\pi}\sigma}\re^{-(x-\mu )^2/(2\sigma^2)}, 
\;\;x\in\bbR .$$ 
It yields an exponential family on $\cX=\bbR$, with a known variance $\sigma^2>0$ 
and a varying mean $\mu$; here the natural parameter is $\tha =\mu$.
Take an exponential WF $\varphi(x)=\re^{\gamma x}$. 
Then 
\beq\label{eq:EthaG}
E_\vphi (\tha)=\dfrac{1}{\sqrt{2\pi}\sigma}\diy\int \re^{\gamma x}\re^{-(x-\tha)^2/(2\sigma^2)}{\rd}x
=\re^{\tha\gamma+\sigma^2\gamma^2/2}\eeq
and the weighted Fisher information is calculated as 
\beq\label{eq:wFIfGaMu}I^\rw_{\varphi}(\tha)=\left(\sigma^{-2}
+\gamma^2\right)\re^{\tha\gamma+\sigma^2\gamma^2/2}.\eeq 
The bound \eqref{wightKRs} takes the form 
\beq\label{eq:wCR4Ew}\beal
\bbE_\tha\Big[\vphi^{(n)} (\bX_1^n)\big(\tha^*(\bX_1^n)-\tha\big)^2\Big]
 \geq\dfrac{\sigma^2}{n}
\dfrac{\left[\re^{n(\tha\gamma+\sigma^2\gamma^2/2)}+b'(\tha )\right]^2}{
\re^{n(\tha\gamma+\sigma^2\gamma^2/2)}(1+n\sigma^2\gam^2)}.\ena\eeq
 
As examples, consider the following estimators for $\theta$:
\beq\beal
\tha^*_1(\bX_1^n)=\dfrac{1}{n}\sum\limits_{i=1}^n X_i,\quad
\tha^*_2(\bX_1^n)=-\sigma^2\gamma+\dfrac{1}{n}\sum\limits_{i=1}^n X_i. %\quad
%{\wh\tha}_1(\bX_1^n)=\dfrac{1}{n E(\tha)}\left(\sum\limits_{i=1}^n X_i\re^{\gamma X_i}\right).
\ena\label{7.41}\eeq
For estimator $\tha^*_1$, the bias $b_1(\tha)=\gam\sigma^2\re^{n(\tha\gamma+\sigma^2\gamma^2/2)}$, and the derivative 
$b'_1(\tha )=n\gam^2\sigma^2\re^{n(\tha\gamma+\sigma^2\gamma^2/2)}$. As a result, 
$\tha^*_1$ yields equality in \eqref{eq:wCR4Ew}:
\beq\label{eq:7.5}\beal \bbE_{\tha}\left[\Big(\tha^*_1(\bX_1^n)
-\tha\Big)^2\exp\Big(\gam\sum\limits_{i=1}^nX_i\Big)\right]\\
\qquad\qquad =
\dfrac{\sigma^2}{n}\re^{n(\tha\gamma+\sigma^2\gamma^2/2)}(1+n\gam^2\sigma^2).
\ena\eeq

For estimator $\tha^*_2$, the bias $b_2(\tha)=0$, and the weighted mean  value is
$$\bbE_{\tha}\left[\tha^*_2(\bX_1^n)
\exp\Big(\gam\sum\limits_{i=1}^nX_i\Big)\right]=\tha\re^{n(\tha\gamma+\sigma^2\gamma^2/2)}.$$ 
The LHS of \eqref{eq:wCR4Ew} has the form
$$\bbE_{\tha}\left[\Big(\tha^*_2(\bX_1^n)
-\tha\Big)^2\exp\Big(\gam\sum\limits_{i=1}^nX_i\Big)\right]=\dfrac{\sigma^2}{n}
\re^{n(\tha\gamma+\sigma^2\gamma^2/2)},$$
whereas the RHS equals \ 
$\dfrac{\sigma^2}{n}\dfrac{\re^{n(\tha\gamma+\sigma^2\gamma^2/2)}}{1+n\gam^2\sigma^2}$.\vskip .5cm 
%which decays as $\dfrac{1}{n^2}$. 
 
Next, $s(\tha)=\re^{\sigma^2\gamma^2/8+\gam\tha /2}$, and estimator $\tha^*_1$ in
\eqref{7.41} yields $c_1(\tha)=\dfrac{1}{2}\sigma^2\gam s(\tha)^n$, $c'_1(\tha )
=\dfrac{1}{4}n\sigma^2\gam^2s(\tha )^n$. Consequently, 
\beq\beal
\bbE_\tha\bigg[\vphi^{(n)} (\bX_1^n)\Big(\tha_1^*(\bX_1^n)-\tha\Big)^2\bigg]
\geq\dfrac{\sigma^2}{n}\re^{n(\gam\tha +\gam^2\sigma^2/4)}
(1+n\gam^2\sigma^2/4)^2.\end{array}\eeq
This bound may be unattainable, and the question of a related estimator with a minimal weighted variance 
remains open. 

{\bf{(B2)}} Consider a family with a scaled parameter: 
$p_\tha (x)=\dfrac{1}{\tha}g\left(\dfrac{x}{\tha}\right)$. Consider $n$ 
IID observations
from PDF $p_\tha$ and assume that the WF $\vphi^{(n)}$ is as in \eqref{eq:vphind}. Then the
weighted Fisher information (cf. \eqref{eq:wFInn}) takes the form
\beq\label{eq:7.7}I^\rw_{\vphi^{(n)}}(\tha)=\dfrac{n}{\tha^2}\diy\int\vphi (x)\left(1
+x\dfrac{g'(x)}{g(x)}\right)^2g(x).\eeq
It is natural to scale the WF with $\tha$, i.e. select
$\varphi(x)=\re^{\gamma x/\tha}$ and $p_\tha (x)=\dfrac{1}{\sqrt{2\pi}\tha}
\re^{-x^2/(2\tha^2)}$. Then $E(\tha)=\re^{\gamma^2/2}$ and
\beq\label{eq:7.8}
I^\rw_{\vphi}(\tha)=\dfrac{1}{\tha^2}\re^{\gam^2/2}(2+4\gamma^2+\gamma^4).\eeq
Substituting \eqref{eq:7.7} into \eqref{wightKRs} yields an explicit form of the weighted
Cram\'er--Rao bound. 
\vskip 1cm

\section{Appendix C. Examples of weighted entropies and divergences} 

In this section we discuss some selected 
examples where we provide simple calculations for weighted 
entropies and divergences. The thrust here is that if we know that PDFs $p_\tha$, 
$p_{\tha'}$ belong to a
specific exponential family then a weighted divergence may be expressed in terms of
two functional parameters: $F^*(\tha )$ (see \eqref{lognor}) and $E(\tha)=\diy\int\vphi\,p_\tha
=\bbE_\tha [\vphi (X)]$ (see \eqref{2.30}--\eqref{2.33} and \eqref{4.36}
--\eqref{eq:4.18}). To avoid repeated calculations, we will focus on weighted Shannon and 
R\'enyi entropies and weighted Kullback-Leibler, Chernoff and Bhattacharyia divergences.
\vskip .5cm %\end{document}

{\bf{(C1)}}  {\it Exponential distributions.} Here $\cX =(0,\infty )$, and we have an
exponential  PDF family $p_\tha =p^\rExp_\lam$, with the natural parameter $\tha =\lam >0$. 
More precisely, the PDF is
$$p^\rExp_\lam (x)=\lam \re^{-\lam x},$$ 
for $x>0$, with $k(x)=0$, $t(x)=-x$ and  the log-normalizer $F(\tha)=-\ln\lam$.
Next, $E_\vphi (\tha)=\lam \wvphi (\lam )$ where $\wvphi$ stands for the Laplace 
transform of WF $\vphi$: $\wvphi(\lam)=\diy\int_0^{\infty} \vphi(x) \re^{-\lam x}{\rd}x$ .
Consequently, $F^*(\tha )=\ln\wvphi(\lam )$; see \eqref{lognor}.
Assume that, for a given $\lam >0$, both $\wvphi(\lam )$ and the derivative 
$\wvphi\,'(\lam )$ exist, and $\wvphi\,'(\lam )=-\diy\int_0^\infty x\vphi (x)\re^{-\lam x}\rd x$.
Then, according to (\ref{2.30}),
\beq\beal
{\tt K}^\rw_{\vphi}(p^\rExp_{\lam}||p^\rExp_{\lam'})=\lam\left(\ln\lam -\ln\lam'\right)
\wvphi(\lam)-\lam (\lam'-\lam)\wvphi'(\lam).\ena\eeq

Similarly, in agreement with \eqref{2.344} and \eqref{eq:REnt1},  
\beq\beac
h^\rw_{\vphi}(p^\rExp_\lam)=-\lam\Big[(\ln\,\lam)  \wvphi(\lam)+\lam\wvphi'(\lam)\Big],\ena\eeq
and
\beq\beac r^\rw_{\vphi ,\alpha}(p^\rExp_\lam )=\dfrac{\lam\wvphi(\lam )}{1-\alpha}
\ln\left[\dfrac{\wvphi(\alpha\lam)}{\wvphi(\lam )}\lam^{\alpha -1}\right].\ena\eeq

Furthermore, in accordance with \eqref{ChernCo}--\eqref{1.3.11},
\beq\beac
{\tt C}^\rw_{\vphi ,\alpha}(p^\rExp_\lam ,p^\rExp_{\lam'})=
-\ln\Big[\lam^\alpha{\lam'}^{1-\alpha}\wvphi\big(\alpha\lam +(1-\alpha )\lam'\Big]
+\ln\big[\lam\wvphi (\lam )\big]\ena\eeq
and
\beq\beac{\tt A}^\rw_\vphi (p^\rExp_\lam,p^\rExp_{\lam'})
=-\ln\Big[\sqrt{\lam\lam'}\,\wvphi\big((\lam +\lam')/2\big)\Big]+\ln [\lam\wvphi (\lam )].\ena\eeq\vskip .5cm
\vskip .5cm

{\bf{(C2)}} {\it Poisson distributions.} This example provides an exponential family of 
probability mass functions $p_\tha =p^\rPo_\lam$. Here  
$$p^\rPo_\lam (l)=\re^{-\lam}\dfrac{\lam^l}{l!}$$ 
where $l=0,1,\ldots$, $\lam >0$, the natural parameter is $\tha =\ln\,\lam$ and 
$t(l)=l$, $k(l)=-\ln\,(l!)$. The log-normalizer is $F(\tha)=\lam$.

In this example we focus on the WF $\vphi (l)=\re^{\gam l}$: then $E_\vphi (\lam)
=\re^{\lam (\re^\gam -1)}$
and
\beq\beal h^\rw_\vphi (p^\rPo_\lam )=E_\vphi (\lam )\Big[\lam-\lam\re^\gam\log \lam
+\re^{-\lam\re^\gam}\sum\limits_{l\geq 0}\dfrac{(\lam\re^\gam)^l}{l!}\ln\,(l!)
\Big].\ena\eeq
Further, 
\beq{\tt K}^\rw_\vphi (p^\rPo_\lam\| p^\rPo_{\lam'})=E_\vphi (\lam )
\Big[\lam'-\lam+(\lam\re^\gam )\log\,\dfrac{\lam}{\lam'}\Big]\eeq
and
\beq\beal h^\rw_{\vphi ,\alpha}(p^\rPo_\lam )=\dfrac{E_\vphi (\lam )}{1-\alpha}
\bigg\{\ln\bigg[\sum\limits_{l\geq 0}
\dfrac{\re^{\gam l}\lam^{\alpha l}}{(l!)^\alpha}\bigg] -\lam\big(\alpha
-\re^\gam +1\big)\bigg\}.\ena\eeq

%Take $\alpha=\dfrac{1}{2}$, then 
%$${\tt U}^\rw_{F,\alpha}(p_\tha,p_{\tha'})=\dfrac{1}{2}(\re^{\tha/2}-\re^{\tha'/2})^2.$$
Next, in view of (\ref{4.36}) we obtain 
\beq\beal
{\tt C}^\rw_{\vphi ,\alpha}(p^\rPo_\lam ,p^\rPo_{\lam'})
=-\re^\gam\lam^\alpha{\lam'}^{1-\alpha}+\alpha\lam+(1-\alpha )\lam'+\lam (\re^\gam -1)
\ena\eeq
and 
\beq\beal{\tt A}^\rw_\vphi (p^\rPo_\lam,p^\rPo_{\lam'})
=-\re^\gam\sqrt{\lam\lam'}+\dfrac{1}{2}(\lam+\lam') +\lam (\re^\gam -1).\ena\eeq\vskip .5cm

{\bf{(C3)}} In the scalar Gaussian case where $X\sim$N$(\mu,\sigma^2)$, we have a family of normal PDFs:
$p_\tha=p^{\rm{Ga}}_{\mu ,\sigma^2}$, with $\cX =\bbR$, $\mu\in\bbR$ and $\sigma^2>0$. Here
$$p^{\rm{Ga}}_{\mu ,\sigma^2}(x)=\dfrac{\exp\big[-(x-\mu )^2/(2\sigma^2)\big]}{
(2\pi\sigma^2)^{1/2}},\quad x\in\bbR.$$
Next, $t(x)=(x,x^2)$, $k(x)=0$, and the natural parameter is $\tha=\big(\mu/\sigma^2,-1/(2\sigma^2)\big)$.
Consequently, the log-normalizer has the form $F(\tha)
=\dfrac{\mu^2}{2\sigma^2}+\dfrac{1}{2}\ln(2\pi\sigma^2)$.

For the weighted Kullback-Leibler divergence, Eqn \eqref{2.30} yields
\beq\label{eq:wKLGa}\beal {\tt K}^\rw_{\vphi}(p^{\rm{Ga}}_{\mu,\sigma^2}||p^{\rm{Ga}}_{\mu',\sigma'^{\,2}})\\
\quad =E_0(\mu ,\sigma^2)\left[\ln\,\dfrac{\sigma'}{\sigma}-\dfrac{1}{2}\left(\dfrac{\mu^2}{\sigma^2}
-\dfrac{\mu'^{\,2}}{\sigma'^{\,2}}\right)\right]\\
\qquad +\left(\dfrac{\mu}{\sigma^2}-\dfrac{\mu'}{{\sigma'}^2}\right)E_1(\mu,\sigma^2)
-\dfrac{1}{2}\left(\dfrac{1}{\sigma^2}-\dfrac{1}{\sigma'^{\,2}}\right)E_2(\mu,\sigma^2).\end{array}\eeq
Here and below, coefficients $E_i (\mu ,\sigma^2)$, $i=0,1,2$, are given by
\beq\label{eq:EsGaS}E_i(\mu ,\sigma^2)=\int_\bbR\vphi (x)p^{\rm{Ga}}_{\mu ,\sigma^2}(x)x^i\rd x,\quad i=0,1,2\eeq
(and $E_0(\mu ,\sigma^2)$ coincides with $E_\vphi (\tha )$.)
\vskip .5cm

According to Eqn \eqref{2.344} and  \eqref{2.33}, the weighted Shannon and R\'enyi 
entropies take the form
\beq\label{eq:wEGa}\beacl h^\rw_{\vphi}(p^{\rm{Ga}}_{\mu ,\sigma^2})
&=\; E_0(\mu ,\sigma^2)\ln\,\Big[(2\pi\sigma^2)^{1/2}\Big]\\
\;&\quad+\dfrac{1}{2}\Big[\mu^2E_0 (\mu,\sigma^2)-2\mu E_1 (\mu,\sigma^2)
+E_2 (\mu,\sigma^2)\Big]\sigma^{-2}.\ena\eeq
%G_2(\mu ,\sigma^2)=\diy\int\vphi (x)x^2g(x)\rd x$$  
and
\beq\label{2.433}\beal r^\rw_{\vphi ,\alpha}(p^{\rm{Ga}}_{\mu ,\sigma^2})
=\dfrac{E_0 (\mu ,\sigma^2)}{1-\alpha}\bigg[\dfrac{1-\alpha}{2}\ln (2\pi\sigma^2)\\
\qquad\qquad-\dfrac{1}{2}\ln\,\alpha+\ln E_0\left(\mu,\sigma^2/\alpha 
\right)-\ln E_0(\mu,\sigma^2)\bigg].\ena\eeq

Further, 
\beq\label{eq:CheGaDiv}\beal
{\tt C}^\rw_{\vphi ,\alpha}(p^\rGa_{\mu ,\sigma^2} ,p^\rGa_{\mu',{\sigma'}^2})
=\alpha \ln\,\sigma +(1-\alpha)\ln\,\sigma'-\ln\osigma\\
\quad +\dfrac{1}{2}\left[\alpha\mu^2\sigma^{-2}
+(1-\alpha){\mu'}^2{\sigma'}^{-2} 
-\omu^2\osigma^{-2}\right]-\ln E_\vphi (\tha )+\ln E_\vphi (\otha).\ena\eeq
Here $\osigma^2=\left(\alpha\sigma^{-2}+(1-\alpha ){\sigma'}^{-2}\right)^{-1}$,
$\omu =\left[\alpha\mu\sigma^{-2}+(1-\alpha )\mu'{\sigma'}^{-2}\right]\osigma^2$, 
and $\otha$ stands for the natural parameter $\left(\omu\,\osigma^{-2},-\osigma^{-2}/2\right)$.
Consequently, $E_\vphi (\otha )=\diy\int\vphi p^\rGa_{\omu ,\osigma^2}$.

Furthermore, 
\beq\label{eq:BhaGaDiv}\beal{\tt A}^\rw_\vphi (p^\rGa_{\mu ,\sigma^2} ,p^\rGa_{\mu',{\sigma'}^2})
=\dfrac{1}{2}\ln\,(\sigma\sigma')-\ln \osigma\\
\qquad +\dfrac{1}{4}\left(\dfrac{\mu^2}{\sigma^2}+\dfrac{{\mu'}^2}{{\sigma'}^2} \right)
-\dfrac{\omu^2}{2\osigma^2}+\ln E_\vphi (\otha)-\ln E_\vphi (\tha)\ena\eeq
where $\osigma^2=\dfrac{2\sigma^2{\sigma'}^2}{\sigma^2+{\sigma'}^2}$, 
$\omu =\dfrac{1}{2}\left(\dfrac{\mu}{\sigma^2}+\dfrac{\mu'}{{\sigma'}^2}\right)\osigma^2$
and $E_\vphi (\otha )=\diy\int\vphi p^\rGa_{\omu ,\osigma^2}$ has the same meaning as above.

%in accordance with \eqref{2.33}. 

In the case of an exponential WF $\vphi (x)=\re^{\gam x}$, $\gam\in\bbR$, then
\beq\label{eq:EGagamS}\beac E_0(\mu ,\sigma^2)=\re^{\mu\gam +\gam^2\sigma^2},\;\;E_1(\mu ,\sigma^2 )
=(\gam\sigma^2 +\mu )E_0(\mu ,\sigma^2),\\
E_2(\mu ,\sigma^2)=\Big[\sigma^2 +(\gam\sigma^2 +\mu)^2\Big] E_0(\mu ,\sigma ^2). \ena\eeq

Consequently, \eqref{eq:wEGa}, \eqref{eq:wKLGa} and \eqref{2.433} are specified as follows
\beq\beal {\tt K}^\rw_{\vphi}(p^{\rm{Ga}}_{\mu,\sigma^2}||p^{\rm{Ga}}_{\mu',\sigma'^2})
=\dfrac{1}{2}\re^{\mu\gam +\gam^2\sigma^2}\Big[\ln\,\dfrac{\sigma^2}{{\sigma'}^2}\\
\qquad -\big(\mu^2\sigma^{-2}-\mu'^{\,2}\sigma'^{\,-2}\big)
+\big(\mu\sigma^{-2}-\mu'\sigma'^{\,-2}\big) (\gam\sigma+\mu )\\
\qquad -\dfrac{1}{2}\big(\sigma^{-2}-\sigma'^{\,-2}\big)\big(\sigma^2 +(\gam\sigma^2 +\mu)^2\big)\Big]
\ena\eeq
\beq\beal h^\rw_\vphi (p^{\rm{Ga}}_{\mu ,\sigma^2})=\dfrac{1}{2}\re^{\mu\gam +\gam^2\sigma^2}\Big(
\ln\Big[(2\pi e)\sigma^2\,\Big] +\gam^2\sigma^2\Big),\ena \eeq
and 
\beq\label{2.57}\bear 
r^\rw_{\vphi ,\alpha}(p^{\rm{Ga}}_{\mu ,\sigma^2})
=\dfrac{1}{2}\re^{\mu\gam +\gam^2\sigma^2}\Big[\ln\big(2\pi\sigma^2\big)
-\dfrac{d\ln\alpha}{1-\alpha}+\dfrac{1}{\alpha}\gam^2\sigma^2\Big].
\ena\eeq

To conclude, for the WF $\vphi (x)=\re^{\gam x}$, 
formulas \eqref{eq:CheGaDiv} and \eqref{eq:BhaGaDiv} admit a further specification,
with $\ln E_\vphi (\otha)=\omu\gam+\osigma^2\gam^2/2$. 
\vskip .5cm

{\bf{(C4)}} For the $d-$dimensional Gaussian case, $X\sim$N$(\mu,\Sigma)$, 
where $\mu=(\mu_1,\ldots ,\mu_d)\in\bbR^d$ and $\Sigma$ is a positive-definite $(d\times d)$ matrix.
The
sample space $\cX=\bbR^d$, and we write $p_\tha =p^\rGa_{\mu ,\Sigma}$
where
$$p^\rGa_{\mu,\Sigma}(x)=\dfrac{\exp\left[-(x-\mu )
\Sigma^{-1}(x-\mu)^{\rT}/2\right]}{((2\pi^d)\rdet\,\Sigma)^{1/2}},\; x=(x_1,\ldots ,x_d)\in\bbR^d.$$
Next, $t(x)=(x,x^\rT x)$, $k(x)=0$, $x\in\bbR^d$, and the natural parameter is
$\tha=\left(\mu\Sigma^{-1},-\Sigma^{-1}/2\right)$. Consequently, the log-normalizer has the form
$$F(\tha)=\dfrac{1}{2}\mu\,\Sigma^{-1}\mu^\rT+\dfrac{1}{2}\ln\Big[(2\pi)^d\rdet\,\Sigma\Big].$$

As in the scalar case, we have that
\beq\beal {\tt K}^\rw_{\vphi}(p^{\rm{Ga}}_{\mu,\Sigma}||p^{\rm{Ga}}_{\mu',\Sigma'})\\
\;=\dfrac{1}{2}E_0(\mu,\Sigma )\Big[\ln\dfrac{\rdet\,\Sigma}{\rdet\,\Sigma'}
-\Big(\mu\Sigma^{-1}\mu^\rT-\mu'\Sigma'^{\,-1}\mu'^{\,\rT} \Big)\Big]\\
\quad+\rtr\;\Big[\big(\mu\Sigma^{-1}-\mu'\Sigma'^{\,-1}\big)^\rT 
E_1(\mu,\Sigma )-\dfrac{1}{2}\big(\Sigma^{-1}-\Sigma'^{\,-1}\big)E_2(\mu,\Sigma)\Big],
\ena\eeq
\beq\label{eq:wSEmGa}\beal h^\rw_\vphi (p^{\rm{Ga}}_{\mu ,\Sigma})=\;E_0(\mu ,\Sigma)
\ln\Big[((2\pi)^d\rdet\,\Sigma )^{1/2}\Big] \\
\qquad +\dfrac{1}{2}\rtr\Big[\big(\mu^{\rT}\mu\;E_0(\mu ,\Sigma )
 -2\mu^\rT E_1(\mu ,\Sigma )+E_2(\mu ,\Sigma\big)\Sigma^{-1}\Big],\ena\eeq
and
\beq\label{2.47}\beacl r^\rw_{\vphi, \alpha}(p^{\rm{Ga}}_{\mu ,\Sigma})
&=\dfrac{E_0 (\mu ,\Sigma )}{1-\alpha}\bigg[\dfrac{1-\alpha}{2}\ln\Big[((2\pi)^d\rdet\,\Sigma )^{1/2}\Big]\\
\;&\quad-\dfrac{1}{2}\ln\,\alpha+\ln E_0\left(\mu,\Sigma/\alpha 
\right)-\ln E_0(\mu,\Sigma)\bigg].\ena\eeq
Coefficients  $E_0(\vphi ,\mu,\Sigma )$ (a scalar), $E_1(\vphi ,\mu,\Sigma )$ 
(a vector) and $E_2(\vphi ,\mu,\Sigma )$ (a matrix) are given by
\beq\beac E_0(\mu ,\Sigma)=\diy\int_{\bbR^d}\vphi (x) p^{\rm{Ga}}_{\mu ,\Sigma}(x)\rd x,\; 
E_1(\mu ,\Sigma)=\diy\int_{\bbR^d}\vphi (x)p^{\rm{Ga}}_{\mu ,\Sigma}(x)x\rd x,\\
E_2(\mu ,\Sigma)=\diy\int_{\bbR^d}\vphi (x)p^{\rm{Ga}}_{\mu ,\Sigma}(x)\,(x^\rT x)\,\rd x.\ena\eeq
(Again, $E_0(\mu ,\Sigma)=E(\tha )$.)

Continuing an analogy with the scalar case, we obtain
\beq\label{eq:CheVeGsM}\beal
{\tt C}^\rw_{\vphi ,\alpha}(p^\rGa_{\mu ,\Sigma} ,p^\rGa_{\mu',\Sigma'})=\alpha\ln\,(\rdet\,\Sigma)
+(1-\alpha )\ln\,(\rdet\,\Sigma')-\ln\,(\rdet\oSigma )\\
+\dfrac{1}{2}\left[\alpha\mu\Sigma^{-1}\mu^\rT
+(1-\alpha){\Sigma'}^{-1}{\mu'}^\rT
-\omu\,\oSigma^{-1}\omu^\rT\right]-\ln E_\vphi (\tha )+\ln E_\vphi (\otha).\ena\eeq
Here $\oSigma =\left(\alpha\Sigma^{-1}+(1-\alpha ){\Sigma'}^{-1}\right)^{-1}$,
$\omu =\left[\alpha\mu\Sigma^{-1}+(1-\alpha )\mu'{\Sigma'}^{-1}\right]\oSigma$, 
and $\otha$ stands for the natural parameter $\left(\omu\,\oSigma^{-1},-\oSigma^{-1}/2\right)$.
Consequently, $E_\vphi (\otha )=\diy\int\vphi p^\rGa_{\omu ,\oSigma^2}$.

Next, 
\beq\label{eq:CheVeGsT}\beal{\tt A}^\rw_\vphi (p^\rGa_{\mu ,\Sigma} ,p^\rGa_{\mu',\Sigma'})
=\dfrac{1}{2}\big[\ln(\rdet\,\Sigma )+\ln(\rdet\,\Sigma')\big]-\ln (\rdet\oSigma)\\
\quad +\dfrac{1}{4}\left[\mu\Sigma^{-1}\mu^\rT+{\Sigma'}^{-1}{\mu'}^\rT
-2\omu\,\oSigma^{-1}\omu^\rT\right]-\ln E_\vphi (\tha )+\ln E_\vphi (\otha)\ena\eeq
where $\oSigma =2\left(\Sigma^{-1}+{\Sigma'}^{-1}\right)^{-1}$ and
$\omu =2\left(\mu\Sigma^{-1}+\mu'{\Sigma'}^{-1}\right)\oSigma$, 
while $\otha$ and $E_\vphi (\otha)$ have the same meaning as above

For the WF $\vphi (x)=\re^{x^\rT\gam}$ where 
$\gam =(\gam_1,\ldots ,\gam_d)\in\bbR^d$, the equations for coefficients $E_i(\mu ,\Sigma )$ become:
\beq\beac E_0(\mu ,\Sigma )=\re^{\mu\gam^\rT +\gam\Sigma\gam^\rT /2},\;\;E_1(\mu ,\Sigma )
=(\gam\Sigma +\mu )E_0(\mu ,\Sigma ),\\
E_2(\mu ,\Sigma )=\Big[\Sigma +(\gam\Sigma +\mu)^\rT(\gam\Sigma +\mu )\Big] E_0(\mu ,\Sigma ). \ena\eeq
Consequently, we obtain the following specifications
\beq\beal {\tt K}^\rw_{\vphi}(p^{\rm{Ga}}_{\mu,\Sigma}||p^{\rm{Ga}}_{\mu',\Sigma'})
=\dfrac{1}{2}\re^{\mu\gam^\rT +\gam\Sigma\gam^\rT /2}\\
\qquad\qquad\times\bigg(\ln\dfrac{\rdet\,\Sigma}{\rdet\,\Sigma'}
-\big(\mu\Sigma^{-1}\mu^\rT-\mu'\Sigma'^{\,-1}\mu'^{\,\rT} \big)\\
\quad\qquad +\rtr\;\Big[\big(\mu\Sigma^{-1}-\mu'\Sigma'^{\,-1}\big)^\rT (\gam\Sigma+\mu )\Big]\\
\qquad\quad-\dfrac{1}{2}\rtr\,\Big[\big(\Sigma^{-1}-\Sigma'^{\,-1}\big)
\big(\Sigma +(\gam\Sigma +\mu)^\rT(\gam\Sigma +\mu )\big)\Big]\bigg),\\
\ena\eeq
\beq\beal h^\rw_\vphi (p^{\rm{Ga}}_{\mu ,\Sigma})=\dfrac{1}{2}\re^{\mu\gam^\rT +\gam\Sigma\gam^\rT /2}
\Big(\ln\Big[(2\pi e)^d\rdet\,\Sigma \,\Big] +\gam\Sigma\gam^\rT\Big),\ena \eeq
and
\beq\beacl
r^\rw_{\vphi ,\alpha}(p^{\rm{Ga}}_{\mu ,\Sigma})
&=\dfrac{1}{2}\re^{\mu\gam^\rT +\gam\Sigma\gam^\rT /2}\\
\;&\quad\times\bigg(\ln\Big[(2\pi)^d\,\rdet\,\Sigma\Big]-\dfrac{d\ln\alpha}{1-\alpha}
+\dfrac{1}{\alpha}\gam\Sigma\gam^{\rT}\bigg).
\label{2.57}\ena\eeq

As before, for $\vphi (x)=\re^{x^\rT\gam}$  formulas \eqref{eq:CheVeGsM} and \eqref{eq:CheVeGsT} 
also admit a further specification. where $\ln E_\vphi (\tha )=\omu\gam^\rT+\gam\Sigma\gam^\rT$.

\vskip .5cm

{\bf{(C5)}} For a Gamma-distribution, $X\sim\;$Gam$(\lam,\beta)$ where $\lam >0$
and $\beta >0$. Here $\cX =(0,\infty )$, and the PDF is
$$p^\rGam_{\lam ,\beta}(x)=\dfrac{\beta^\lam}{\Gam (\lam )}x^{\lam -1}\re^{-\beta x},$$
for $x>0$.
Gamma-PDFs form an exponential PDF familiy, with $t(x)=(x, \ln x)$ and $k(x)=0$,
where the natural parameter is $\tha=(-\beta,\lam-1)$. Consequently, the log-normalizer is
$F(\tha )=\ln\dfrac{\Gamma(\lam)}{\beta^\lam}$; 
cf. \cite{Mor}. As above, we use the notation $E_\vphi (\tha )$ for the mean 
$\diy\int_0^\infty\vphi (x)p^\rGam_{\lam ,\beta}(x)\rd x$ of a WF $\vphi$. 
Also set: 
$$L(\lam ,\beta)
=\diy\int_0^\infty\vphi(x)\big[\beta x+(1-\lam)\ln x\big]p^\rGam_{\lam ,\beta}(x)\rd x.$$
 
Then
\beq\beal{\tt K}^\rw_\vphi (p^\rGam_{\lam ,\beta}\|p^\rGam_{\lam',\beta'}) 
=E_\vphi (\tha )\ln\dfrac{\beta^\lam\Gam (\lam ')}{{\beta'}^{\lam'}\Gam (\lam)}
+L(1-\lam +\lam',\beta -\beta')\ena\eeq
\beq\beacl
h^\rw_{\vphi}(p^\rGam_{\lam ,\beta})&= E_\vphi (\tha)\ln\dfrac{\Gamma(\lam)}{\beta^\lam}+L(\lam ,\beta )
\ena\eeq
%in accordance with \eqref{2.344}. 
and %the weighted R\'enyi entropy as 
\beq\beal
r^\rw_{\vphi ,\alpha}(p^\rGam_{\lam ,\beta})
=\dfrac{E_\vphi (\tha )}{1-\alpha}\bigg[\ln\dfrac{E_\vphi (\alpha\tha)}{E_\vphi (\tha )}
+\ln\dfrac{\Gamma(\alpha\lam)}{(\alpha\beta)^{\alpha\lam}} 
-\alpha\ln\dfrac{\Gamma(\lam)}{\beta^\lam}\bigg],\ena\eeq
confirming \eqref{2.30}, \eqref{2.344} and \eqref{2.33}, respectively.

Finally, set $\olam=\olam (\alpha )=\alpha \lambda +(1-\alpha) \lambda'$ and 
$\obeta =\obeta (\alpha )=\alpha\beta+(1-\alpha)\beta'$.
Consider the Gamma-PDF $p^\rGam_{\olam ,\obeta}$ with
the value of the natural parameter $\otha=(-\obeta ,\olam -1)$.
Then
\beq\beal {\tt C}^\rw_{\vphi ,\alpha}\big(p^\rGam_{\lam ,\beta},p^\rGam_{\lam',\beta'}\big)\\
\quad =\ln E_\vphi (\tha )-\ln E_\vphi (\otha)+\alpha\ln\dfrac{\Gamma(\lam)}{\beta^\lam}
+(1-\alpha)\ln\dfrac{\Gamma(\lam')}{{\beta'}^{\lam'}}
-\ln\dfrac{\Gamma(\olam)}{\obeta^\olam}.\ena\eeq
Furthermore, for $\alpha=1/2$ we obtain $\otha=(-(\beta +\beta')/2,(\lam +\lam')/2 -1)$ 
\beq\beal{\tt A}^\rw_{\vphi }\big(p^\rGam_{\lam ,\beta},p^\rGam_{\lam',\beta'}\big)\\
\quad =\ln E_\vphi (\tha )-\ln E_\vphi (\otha)+\dfrac{1}{2}\left[\ln\dfrac{\Gamma(\lam)}{\beta^\lam}
+\ln\dfrac{\Gamma(\lam')}{{\beta'}^{\lam'}}\right]
-\ln\dfrac{\Gamma(\olam)}{\obeta^\olam}.\ena\eeq

\vskip .5cm

\end{document}